\documentclass[12pt]{article}
\usepackage{amsmath,amsfonts,amssymb,theorem}
  \usepackage[matrix, arrow]{xy}

\DeclareMathOperator{\Supp}{Supp}

\DeclareMathOperator{\thr}{th}
\DeclareMathOperator{\LCS}{LCS}

\DeclareMathOperator{\codim}{codim}

\DeclareMathOperator{\mld}{mld}

\DeclareMathOperator{\mult}{mult}

\DeclareMathOperator{\lct}{lct}
\DeclareMathOperator{\LMMP}{LMMP}
\DeclareMathOperator{\ACC}{ACC}

\DeclareMathOperator{\Pic}{Pic}
\DeclareMathOperator{\act}{act}
\DeclareMathOperator{\ldis}{ldis}

 \newcommand{\iso}{\cong}
 
 \newcommand{\pt}{\mathrm{pt.}}

 \newcommand{\rddown}[1]{\left\lfloor{#1}\right\rfloor} 

 \newcommand{\ep}{\varepsilon}

 \numberwithin{equation}{subsection}
 \numberwithin{footnote}{section}

 \newtheorem{mtheorem}[subsection]{Main Theorem}

 \newtheorem{cor}[subsection]{Corollary}
 \newtheorem{lem}[subsection]{Lemma}
 \newtheorem{prop}[subsection]{Proposition}
 \newtheorem{thm}[subsection]{Theorem}
 \newtheorem{conj}[subsection]{Conjecture}
 \newtheorem{adden}[subsection]{Addendum}
 
 \newtheorem{mainp}[subsection]{Main Proposition}

 {
 \theorembodyfont{\rmfamily}
 \newtheorem{defn}[subsection]{Definition}
 \newtheorem{exa}[subsection]{Example}
 \newtheorem{exa-cr}[subsection]{Example--Construction}

 \newtheorem{rem}[subsection]{Remark}
 
 \newtheorem{warn}[subsection]{Caution}
  \newtheorem{constr}[subsection]{Construction}}
  \newtheorem{rem-defn}[subsection]{Definition}

 \newcommand{\qed}{\ifhmode\unskip\nobreak\fi\quad\ensuremath\square}
 \newenvironment{proof}{\paragraph{Proof}}{\par\medskip}

 \newcommand{\ke}[1]{$\acute{\mbox{e}}$}
 \newcommand{\ku}[1]{$\acute{\mbox{u}$}}
 \newcommand{\kl}[1]{$\acute{\mbox{l}}$}
 \newcommand{\kh}[1]{$\acute{\mbox{h}}$}
 \newcommand{\kr}[1]{$\acute{\mbox{r}}$}
 \newcommand{\kx}[1]{$\acute{\mbox{x}}$}
 \newcommand{\ki}[1]{${\^\i}$}

\title{Mld's vs thresholds and flips} 

\author{C. Birkar\thanks{Supported by the EPSRC.}~~~~~~~~~~~~~~~~~~V.V. Shokurov\thanks{Partially supported
by NSF grant DMS-0400832.}
}

\begin{document}
\maketitle

\begin{abstract}
 Minimal log discrepancies (mld's) are related not only to termination of log flips
[\ref{mld's}], and
thus to the existence of log flips [\ref{F}]
but also to the ascending chain condition (acc) of some global invariants and
invariants of singularities in the Log Minimal Model Program (LMMP). In this paper,
we draw clear links between several central conjectures in the LMMP. More precisely,
our main result states that the LMMP, the acc conjecture for mld's and the boundedness
of canonical Mori-Fano varieties in dimension $\le d$
imply the following:
the acc conjecture for $a$-lc thresholds, in particular,
for canonical and log canonical (lc) thresholds in dimension $\le d$;
the acc conjecture for lc thresholds in dimension $\le d+1$;
termination of log flips in dimension $\le d+1$ for effective pairs; and
existence of pl flips in dimension $\le d+2$.
This also gives new proofs of some well-known and new results
in the field in low dimensions:
the acc conjecture holds for $a$-lc thresholds of surfaces;
the acc conjecture holds for lc thresholds of $3$-folds;
termination of $3$-fold log flips holds for effective pairs; and
the existence of $4$-fold pl flips holds.

\end{abstract}

\tableofcontents

  \section{Introduction}

Two main open problems in the Log Minimal Model Program (LMMP) are:
existence and termination of log flips.
Essentially, the latter one is the only problem:
LMMP in dimension $d$, or inductively just termination in dimension $d$
implies existence of log flips in dimension $d+1$ [\ref{F}].
Thus, if we establish termination of log flips in
dimension $d+1$, LMMP will be completed.
On the other hand, the termination follows from
two local (even formal) problems [\ref{mld's}]:
the ascending chain condition ($\ACC$) conjecture for
 minimal log discrepancies (mld's; see
Conjecture~\ref{acc} below), and
their semicontinuity conjecture due to
Florin Ambro
[\ref{Ambro-mld2}, Conjecture~2.4].
Recently, the first author [\ref{B}] reduced a weaker termination in
dimension $d+1$ (e.g., when the log Kodaira dimension is nonnegative),
in particular, termination of log flips in the relative birational case,
to $\ACC$ conjecture for log canonical (lc) thresholds in dimension $d+1$
(see Conjecture~\ref{ACC}) which in its turn
follows from V.~Alexeev's, and brothers'
A. and L. Borisov conjecture.
This implies a weaker version, in particular,
birational one, of LMMP in dimension $d+1$.
In this paper, we establish that the first of
these conjectures, the  $\ACC$ conjecture for
mld's and a rather weak form of V.~Alexeev's, and brothers'
A. and L. Borisov conjecture (see Conjecture~\ref{weak-BAB})
in dimension $d$ imply the weak version of LMMP in dimension $d+1$.
We hope that this version could be useful
to resolve the two above local conjectures about mld's.
We also show that the same conjectures
are naturally related to some other similar
problems in the field.

We use the terminology of [\ref{log models}] [\ref{pl
flips}] [\ref{Isk-Sh}][\ref{K$^+$}];
see also Notation and terminology below.
However we need certain modifications or
generalizations of some well-known notions and
conjectures.

\begin{defn}[cf. [\ref{Isk-Sh}, Definition~1.6(v)] A proper contraction $X\to Z$ of normal
varieties is called a {\em Mori-Fano fibration\/} if
the following conditions hold:

a) $\dim Z< \dim X$;

b) $X$ has only $\mathbb{Q}$-factorial log canonical (lc) singularities;

c) $\rho(X/Z):=\rho(X)-\rho(Z)=1$, where $\rho(\ )$ is
the Picard number; and

d) $-K$ is ample on $X/Z$.

If $Z=\pt$
is a point, $X$ is called a {\em Mori-Fano variety\/}.
We say that $X$ is a {\em canonical\/} Mori-Fano variety if
$X$ has only canonical (cn) singularities.

\end{defn}

Note that by the Kleiman projectivity criterion,
any Mori-Fano fibration and variety are projective.

\begin{conj}[Weak BAB]\label{weak-BAB} The canonical $d$-dimensional
Mori-Fano varieties are {\em bounded\/}, that is,
a coarse moduli space of such varieties is well-defined
and of finite type. 

\end{conj}

BAB abbreviates V.~Alexeev, and brothers
A. and L. Borisov. The conjecture is a very special case of
their conjecture (see [\ref{MP}]). Conjecture \ref{weak-BAB} is established in dimension $\le 3$ in characteristic zero [KMMT] (the case $d=2$ is
classical). Actually, we need a much weaker version of this conjecture, namely, the boundedness of
canonical $d$-dimensional Mori-Fano varieties $X$ such that $K+B\equiv 0$ for some boundary $B\in \Gamma$ where $\Gamma$ is a fixed set
of boundary multiplicities
satisfying the descending chain condition (dcc).

\begin{conj}[ACC for mld's]\label{acc}
Suppose that $\Gamma\subseteq [0,1]$ satisfies the dcc.
Then the following is expected:

{\rm (ACC)~~~~} The following subset of real numbers $\mathbb{R}$
\[ \big\{\mld(P,X,B)\big\vert ~{\mbox{$(X,B)$
is lc, $\dim X= d$, $P \in X$, and $B\in \Gamma$}}\big\}\]
satisfies the ascending chain condition (acc)

A point $P$ can be nonclosed.
Equivalently, we can consider only closed points $P\in X$, and
assume that $\dim X\le d$.
\end{conj}

This conjecture is established in dimension $d\le 2$ [\ref{A}] [\ref{acc-cod2}],
and for some special cases in higher dimensions [\ref{Bo}] [\ref{log models}][\ref{Ambro-mld}].

\begin{defn}[$a$-lc thresholds]\label{a-threshold} Let $a\geq 0$ be a real number, $(X,B)$ be a log pair,
and $H$ be an $\mathbb{R}$-Cartier divisor on $X$. Then the real number or $+/-\infty$:
$$
t=\thr_a(M,X,B)=\sup \{\lambda\in R\mid (X,B+\lambda H)~~~~~~~~~~~~~~~~~~~~~~~~~~~
$$
$$~~~~~~~~~~~~~~~~~~~~~~~~~~~~~~~~~~~~~~
\mbox{ is $a$-lc in codimension $\ge 2$~ [\ref{log models}, 1.3]}\}
$$
is called the $a$-{\em lc threshold of $H$ with respect to\/} $(X,B)$.
In particular, if $a=0$ or $a=1$, the $a$-lc threshold is
the {\em lc threshold\/} or
{\em cn threshold\/} respectively. $\mathbb{Q}$-{\em factorial} threshold means
that we consider only $\mathbb{Q}$-factorial varieties. Similarly, we get the $a$-lc threshold {\em at
a point\/} $P$ (possibly not closed) if
the $a$-lc condition in codimension $\ge 2$ is replaced by
the $a$-lc condition in $P$.

\end{defn}

\begin{rem}\label{th-rem}
Note that if $+\infty>\ldis(X,D)\ge a$,  
and $H> 0$, then
$t\ge 0$, $\sup=\max$, and is a nonnegative
real number (that is, not $+\infty$, cf. [\ref{Isk-Sh}, Remark~1.4,(ii)]).
In this situation, either $(X,D+tH)$ is precisely $a$-lc in codimension 2, that is $\ldis(X,D+tH)=a$, or
$(X,D+tH)$ has reduced components. Behavior of thresholds in codimension $1$
(at divisorial points) is easy. However, when we consider thresholds at a point,  the situation is more complicated (see Example~\ref{behavior} or cf.
the proof of Proposition~\ref{mld-th}).

Note that we need only $a\le 1$ if $\dim X\ge 2$.
Indeed, $\ldis (X,D)\le 1$ always when $\dim X\ge 2$, and
$\ldis (X,D)=+\infty$ when $\dim X \le 1$ because it corresponds
to the empty set (see Notation and terminology:
mld in codimension $\ge 2$. Cf. thresholds at a point in
Definition~\ref{a-threshold} above).
In contrast, for $\ldis (X,B)<a$ and $H>0$, $t<0$ holds.
Usually in applications, $\dim X\ge 2$, $1\ge \ldis (X,B)\ge a$, and
$D=B$ is a boundary, e.g., $D=0$ [\ref{Isk-Sh}][\ref{K$^+$}][\ref{Cheltsov-Park}][\ref{Ein-Mustata}].

\end{rem}

\begin{exa}\label{behavior}
The $a$-lc threshold at $P$ may not be attained at $P$ nor on the boundary. For example: take three planes $S_1,S_2,S_3$
in the space $\mathbb P^3$ passing through a line $L$. Take a closed point $P\in L$ and define $B=\frac{2}{3}S_1+\frac{2}{3}S_2+\frac{2}{3}S_3$.
$L$ is a lc centre for $(\mathbb{P}^3,B)$ but easy computations show that $\mld(P,\mathbb{P}^3,B)=1$. On the other hand,  $(b_1,b_2,b_3)\neq (1,1,1)$.

\end{exa}

So, in general, for the $a$-lc threshold at a point $P$,
either we get a lc centre passing through $P$ or the mld $a$ is attained at $P$.

\begin{conj}[ACC for $a$-lc thresholds]\label{ACC}
Suppose $d\geq 2$ is a natural number, $a\geq 0$, $\{0\}\subseteq \Gamma\subseteq [0,1]$ satisfies the dcc
and $\{0\}\subsetneq S \subset \mathbb{R}$ is a finite (even dcc)  set of nonnegative numbers. Then the
following is expected:

{\rm (ACC)~~~~~}
The subset
$$
\mathcal T_{a,d}(\Gamma,S)=\{\thr_a(M,X,B)|\mbox{ $(X,B)$ is $a$-lc in codimension $\ge 2$, $\dim X= d$},$$
$$
\mbox{~~~~~~~~~~~~~~~~~~~~~~~~~$B\in \Gamma$, $M$ is an $\mathbb{R}$-Cartier
divisor on $X$, and $M\in S$}\}.
$$
of $\mathbb{R}^+\cup \{+\infty\}$ satisfies  the acc; $+\infty$ corresponds to the case $M=0$.

It is also expected that ACC holds for $a$-lc thresholds
at a point $P$, that is, for the set with
the $a$-lc in $P\in X$ (see Definition~1.4).
The latter set is larger. Thus ACC for thresholds
at a point implies ACC for thresholds on a variety, and
in what follows, ACC for thresholds
means at a point.

\end{conj}

 The case when $d=1$ is obvious: the set of
$\mld(P,X,B)=1-\mult_P B$, for prime divisors $P$ on $X$,
satisfies the acc if and only if the multiplicities of
possible $B$ satisfy the dcc.
Similarly, ACC for $a$-lc thresholds at prime divisors
can be easily verified (cf. Example~\ref{d=1} below).

\begin{mtheorem}\label{result}
$\LMMP$, $\ACC$ for mld's and Conjecture \ref{weak-BAB} in dimension $\leq d$
imply the following:
 
{\rm (i)}
$\ACC$ for $a$-lc thresholds in dimension $\le d$;

{\rm (ii)}
$\ACC$ for lc thresholds in dimension $\le d+1$;

{\rm (iii)}
termination of log flips in dimension $\le d+1$ for effective pairs; and

{\rm (iv)}
existence of pl flips in dimension $\le d+2$.

\end{mtheorem}

See the proof in Section~5.
For generalizations of statements (i-ii) in the Main Theorem and of the
corollaries below, see Section~2.

\begin{cor} ACC for  mld's and Conjecture \ref{weak-BAB} for $4$-folds
imply:

{\rm (i)}
$\ACC$ for $a$-lc thresholds in dimension $ 4$;

{\rm (ii)}
$\ACC$ for lc thresholds in dimension $ 5$;

{\rm (iii)}
termination of log flips in dimension $ 5$ for effective pairs; and

{\rm (iv)}
existence of pl flips in dimension $ 6$.

\end{cor}

\begin{proof}
ACC for mld's for $4$-folds implies termination of
$4$-fold log flips [\ref{mld's}, Corollary 5] and thus LMMP for
$4$-folds [\ref{pl flips}, Corollary~1.8].  $\Box$
\end{proof}

\begin{cor}
$\ACC$ for mld's of $3$-folds implies:

{\rm (i)}
$\ACC$ for $a$-lc thresholds in dimension $ 3$;

{\rm (ii)}
$\ACC$ for lc thresholds in dimension $4$;

{\rm (iii)}
termination of 4-fold log flips for effective pairs;
and

{\rm (iv)}
existence of pl flips in dimension $ 5$.
\end{cor}

\begin{proof}
Immediate by Theorem \ref{result}, and
by [\ref{KMMT}], [\ref{3-fold term}] for $\mathbb{Q}$-boundaries and [\ref{log models}] in general.
$\Box$
\end{proof}

   $\ACC$ for mld's of  algebraic surfaces gives a new
proof of the following well known and new results.

\begin{cor}\label{cor-surfaces}
The following hold:

{\rm (i)}
$\ACC$ for $a$-lc thresholds of surfaces;

{\rm (ii)}
$\ACC$ for lc thresholds of $3$-folds;

{\rm (iii)}
termination of $3$-fold log flips for effective pairs; and

{\rm (iv)}
the existence of $4$-fold pl flips.
\end{cor}

\begin{proof}
Immediate by Theorem \ref{result}, and [\ref{A}][\ref{acc-cod2}] and
[\ref{A1}].
$\Box$
\end{proof}

Note that ACC for mld's of surfaces in [\ref{acc-cod2}] is established for
$\mathbb{R}$-boundaries and without using classification.
Thus, for the first time, termination in  (iii) is proved without
classification (cf. [\ref{3-fold term}] [19, proof of 5.1.3 for 3-folds]). However, this termination for $3$-folds
is still partial.

Cn thresholds and, in particular, their ACC is crucial for
the Sarkisov program [\ref{Corti}] [\ref{Matsuki}].
Another similar important invariant , the Sarkisov degree or its inverse
-- the anticanonical threshold  -- can be included into more general ones:
Fano indices (see Cor.~\ref{anticanonical} below) and boundary
multiplicities of log  pairs
for $S_d(global)$ (see Def.~\ref{type} (v) and Weak finiteness~\ref{Weak
finiteness}).
These invariants and results about them are important in
the proof of our Main Theorem and will be discussed in
Sections~2--4. Here we give a sample.

\begin{cor} Let $\Gamma\subset [0,1]$ be a dcc set. Then, there is a finite subset
$\Gamma_f\subset\Gamma$ such that  $\mathcal S_3^0(\Gamma, global)=\mathcal S_3^0(\Gamma_f,
global)$ (see Definition \ref{type}). \end{cor}

In other words, the set of boundary multiplicities which occur
on the following log pairs is finite: the $3$-fold projective log pairs
$(X,B)$ with $B\in \Gamma$,  $K+B\equiv 0$, $(X,B)$ is lc but not klt.

\begin{proof} Immediate by Theorem \ref{main-th} (vi). $\Box$
\end{proof}

\subsection*{Notation and terminology}

In this paper, a log pair $(X/Z,B)$ consists of normal algebraic varieties
$X,Z$
over a base field $k$ of characteristic $0$, e.g, $k=\mathbb{C}$, where $X/Z$
is a projective morphism, and an $\mathbb{R}$-boundary $B$ (i.e., a divisor
with multiplicities in
$[0,1]$) such that $K+B$ is $\mathbb{R}$-Cartier. Of course, some results
hold or are expected over any field, e.g., $\ACC$ for
$a$-lc thresholds holds in Corollary~\ref{cor-surfaces} (i). We consider
the log minimal model program
($\LMMP$)[\ref{log models}, 5.1] in dimension $d$ in the category of lc pairs of dimension $d$.

An \emph{effective} log pair is a log pair $(X/Z,B)$ [\ref{B}] such that  $K+B\equiv M/Z$ for some $\mathbb{R}$-divisor $M\ge 0$. This property is preserved under any log flip or divisorial contraction. A variety $X$ is of {\em{ Fano type}} (FT) if there an $\mathbb{R}$-boundary $B$ such
that $(X,B)$ is a klt weak log Fano.

A property holds at a point $P\in X$ means that that property holds at the point $P$
but not necessarily in a neighbourhood of $P$. On the other hand, a
property holds near $P$ means that that property holds in an open neighbourhood
of $P$.

If $(X,B)$ is lc, then
$$
1-\mld(P,X,B)=\max\{\mbox{$\mult E$ in $B_W$}|\mbox{$E$ is a prime divisor}~~~~~
$$
$$
~~~~~~~~~~~~~~~~~~~~~~~~~~~~~~~~~~~~~~~~~~~~~~~~~~~~~~~~\mbox{on $W$ and} ~ f(E)=P\}
$$
for any (crepant) log resolution $W\to X$  where $K_W+B_W=f^*(K+B)$ and $\mult$
stands for the multiplicity function on divisors. We define
$$
\ldis(X,B)=\min \{\mld(P,X,B)\mid \mbox{$P\in X$ is of codimension $\ge 2$}\}
$$
 
We say $(X,B)$ is $a$-lc at $P\in X$ if $\mld(P,X,B)\ge a$.
This implies, in particular, that $(X,B)$ is lc near $P$ [\ref{log
models}, Corollary 1.5].

For a set  $\Gamma\subset \mathbb{R}$ and an $\mathbb{R}$-divisor $D$ on a variety, by $D\in\Gamma$ we
mean that the (nonzero) multiplicities of $D$ are in $\Gamma$.


\section{Acc of mld's and thresholds}

For $\mathbb{R}$-divisors on $X$, we have the well-known {\em order\/}:

$D_1\ge D_2$ if $D_1-D_2\ge0$, that is effective.

On the other hand, the topology and the following natural norm,
the {\em maximal absolute value norm\/}, are well-known:
if $D=\sum d_i D_i$, where $d_i\in \mathbb{R}$, and $D_i$ are
distinct prime divisors on $X$, set
$$
\parallel D\parallel=\max \parallel d_i\parallel.
$$

In particular, limits of {\em divisors\/}
are limits in the norm.

\begin{mainp}\label{transform}
We assume $\ACC$ for mld's in dimension $d$.
Let  $\Gamma\subset [0,1]$ be a dcc subset, and
$a$ be a positive real number.
Then, there exists
a real number $\tau>0$ (depending also on $d$)
satisfying the following {\em upper approximation property\/}:
if  $(X,B)$ and $(X,B')$ are two log pairs with a point
$P\in X$ (not necessarily closed) such that

(1)
$\dim X=d$;

(2)
$B\le B'$, $\parallel B-B'\parallel <\tau$, $B'\in \Gamma$;
and

(3) $\mld(P,X,B)\ge a$, that is, $(X,B)$ is $a$-lc at $P$, and
$K+B'$ is $\mathbb{R}$-Cartier; we can omit the last assumption
when $X$ is $\mathbb{Q}$-factorial; and

(4)
$(X,B')$ is lc in a neighborhood of P;

then
$\mld(P,X,B')\ge a$ and $(X,B')$ is also $a$-lc at $P$.
\end{mainp}

Note that by [\ref{log models}, Lemma~1.4] the assumption (4) of the proposition is
equivalent to the lc property of $(X,B')$ at $P$, that is,
to the inequality $\mld(P,X,B')\ge 0$.
To prove the proposition we need the following
general fact.

\begin{lem}[Continuity]\label{continuity}
Suppose that the pairs $(X,B)$ and $(X,B')$ are lc in a neighborhood of a
point $P$.
Then, $a'=\mld(P,X,B')$ and  $a=\mld(P,X,B)$ are real numbers $\ge 0$,
and,
for any real number $x$  in the interval $[a',a]$
there exist two real numbers $\alpha,\beta\ge 0$ such that
$\alpha+\beta=1$, and
$\mld(P,X,\alpha B+\beta B')=x$.
\end{lem}

\begin{proof}
Let $D=B'-B$, and $a\ge a'$.

By the lc property, both mld's are real numbers $\ge 0$.
Then, the last statement holds for the pairs $(X,B+tD)$
with $t=0$ and $t=1$ for which
respectively $B+tD=B$ and $B'$.
By the convexity in [\ref{log flips}, (1.3.2)] the same holds for any
$t$ in the interval $[0,1]$: $(X,B+t D)$ is lc
near $P$, and
$\mld(P,X,B+t D)$ is a real number $\ge 0$.

Let $s=\sup\{t\in [0,1]|\mld(P,X,B+tD)\ge x\}$.
The set of such $s$ is not $\emptyset$
because $a=\mld(P,X,B)\ge x\ge a'$.
We claim that $\mld(P,X,B+sD)=x$.
Put  $\mld(P,X,B+sD)=y$.
Since the last mld is a real number (not $-\infty$)
there is a log resolution $f\colon W\to X$ of $(X,B+sD)$  on which
the mld is attained on divisors:
$$
1-y=\max\{\mbox{$\mult E$ in $B_W+sf^*D$}|\mbox{$E$ is a prime divisor on
$W$ and} ~ f(E)=P\},
$$
where $B_W$ denotes the crepant pull-back of $B$, that is,
$K_W+B_W=f^*(K+B)$.
Then for any $t$,
$\mld (P,X,B+tD)\le 1-m(t)$ with
$$
m(t):=\max\{\mbox{$\mult E$ in $B_W+tf^*D$}|\mbox{$E$ is a prime divisor
on $W$ and} ~ f(E)=P\}
$$
because to calculate the mld one needs to consider the inf (of log discrepancies) for all resolutions.
Note that  $m(t)$ is a piecewise linear and continuous real-valued function
of $t$. This follows from the linear property of $\mult _E(B_W+tf^*D)$
with respect to $t$.

If $y<x$, then for any $t$ sufficiently close to $s$, and
in particular, for such $t<s$,
$\mld (P,X,B+tD)<x$ too, which contradicts our constructions
when $s>0$.
Thus $s=0$ or $y\ge x$, and actually $y=x$ in both cases
by the following stability property:
{\em if $\mld(P,X,B+sD)>x$, and $(P,X,B+tD)$ is lc in
a small neighborhood of $s$ in $[0,1]$, then $\mld(P,X,B+tD)>x$
too in some small neighborhood of $s$ in $[0,1]$\/}.
Note that $s<1$ for $y>x$.
After taking a log resolution:
the stability follows from its log nonsingular version:
for a normal crossing $\mathbb{R}$-subboundary $C$ on $X$ and any point $P$,
if $\mld (P,X,C)>x\ge 0$, then the same holds for any
small perturbation of $C$ in the nonreduced part,
that is, the perturbation of only multiplicities of $C<1$. In fact let $C=\rddown{C}+\sum c_iD_i$ and we may assume that all the components pass through $P$. Then, $\mld(P,X,C)=\codim P-\mu_P\rddown{C}-\sum c_i>x$ where
  $\mu_P\rddown{C}\in \mathbb N$ is the multiplicity of the reduced part $\rddown{C}$ at $P$. Thus, $\mld(P,X,C+\sum \ep_iD_i)=\codim P-\mu_P \rddown{C}-\sum (c_i+\ep_i)>x$ where $|\ep_i|$ are small enough.

Now by taking $\beta=s$ and $\alpha=1-\beta$, we get
$$x=\mld(P,X,(\alpha+\beta)B+\beta(B'-B))=\mld(P,X,\alpha B+\beta B').$$

$\Box$
\end{proof}

\begin{proof}{\bf{of Proposition \ref{transform}}}

Suppose that the proposition does not hold.
Then, there exists a sequence of positive real numbers
$\tau_{1}>\tau_{2}>\dots$ with $\lim_{i\to +\infty} \tau_i=0$, and
a sequence of $d$-dimensional log pairs $(X_i,B_i), i=1,2,\dots,$
such that the proposition does not hold for $\tau_i$ on
$(X_i,B_i)$ in a point $P_i\in X_i$.
In other words, there exists
$B_i'\in \Gamma$ on $X_i$ under (2-4) with $\tau=\tau_i$, and

$$\mld(P_i,X_i,B_i)\geq a ~~~ {\rm but}~~~ \mld(P_i,X_i,B_i')<a$$

We now construct a new sequence of $d$-dimensional log pairs $(T_i,A_i)$ and points $Q_i\in T_i$
such that $a_i=\ldis(T_i,A_i)$ is strictly increasing with $i$ and such that $\Omega$, the set of
multiplicities of all boundaries $A_i$, satisfies the dcc.

By $\ACC$ for mld's the set $\{a_i'=\mld(P_i,X_i,B_i')\}$ has a maximum which is less than $a$. We can assume
that this maximum is equal to $\mld(P_1,X_1,B_1')$. Put $(T_1,A_1):=(X_1,B_1')$ and $Q_1=P_1$
and let $a_1=\mld(Q_1,T_1,A_1)$. Note that $a_1$ is a real number $\ge 0$ and
$a_1=-\infty$ is
impossible by (4).

Suppose that we have already constructed $(T_j,A_j)$ for $1\leq j\leq i-1$. Since $\Gamma$ satisfies the dcc, we can
choose $\tau_k$ such  that there are no multiplicities of $A_j$, for $1\leq j\leq i-1$, in
$(r-\tau_k,r)$ for any $r\in \Gamma$. Take $(T_i,A_i):=(X_k,\alpha B_k+\beta B_k')$, for some
$\alpha, \beta>0$ with $\alpha+\beta=1$, and $Q_i=P_k$ such that
$$ \frac{a_{i-1}+a}{2}< a_i=\mld(Q_i,T_i,\Delta_i)<a.$$

The existence follows from Lemma \ref{continuity} with
$X=X_k, B=B_k, B'=B_k'$, $a'=a_k'$, and any $x$
in the interval $(a',a)$
(here the $a$ in the proposition!) but
$a=a_k$ in the lemma (not the $a$ in the proposition).
Such $x$ exists because, by construction and
assumptions (2-4),
$a_{i-1},a_k'<a$ (for both $a$).

Also by construction for every real number $\ep>0$,
almost all (expect for finitely many) multiplicities
of $\Omega$ belong to intervals
 $(r-\ep,r]$ where $r\in\Gamma$.
This implies that $\Omega$ satisfies the dcc because $\Gamma$ does so. On
the other hand, the set of mld's $\{a_i\}$ does not satisfy the acc which contradicts the
$\ACC$ for mld's. $\Box$ \end{proof}

\begin{prop}\label{aa}
Under the assumptions of Proposition \ref{transform},
let $Y\to X$ be an extremal divisorial extraction
such
that

(1)
$K_Y+B_Y'+(1-a)E$ is $\mathbb{R}$-Cartier,

where $E$ is
the exceptional reduced divisor (but
not necessarily irreducible), and
$B_Y'$ is the birational transform of $B'$ on $X$.
Then $K_Y+B_Y'+(1-a)E$ is seminegative$/X$.
\end{prop}

\begin{exa}\label{bb} The typical situation where we apply
the proposition is as follows.
Let $(X_i,B_i), i=1,2, \dots,$ be a sequence of $d$-dimensional
plt log pairs such that

a) $a_1=\ldis(X_1,B_1)\ge \dots\ge a_i=\ldis(X_i,B_i)\ge\dots>0$
with

b) $a=\lim_{i\to\infty}a_i>0$; and

c) $B_1\le \dots\le B_i\le \dots $ with

d) $B=\lim_{i\to \infty} B_i \in R$.

The c)-d) means that there exist prime divisors
$D_{i,k}, k=1,\dots,n$, on each $X_i$ such that
every
$$
B_i=\sum_{k=1}^{n} b_{i,k}D_{i,k};
$$
and, for every $k=1,\dots,n$,

c') $b_{1,k}\le \dots \le b_{i,k}\le \dots $ with

d') $b_k=\lim_{i\to \infty} b_{i,k} \in R$

(cf. types in Definition \ref{type} below).

In particular, $B=\sum b_kD_{i,k}$ approximates
$B_i$  on $X_i$,
that is, for any $\tau>0$, $\parallel B_i-B\parallel<\tau$
for all $i\gg 0$ (divisors $B$ on $X_i$ have the same
type $(b_1,\dots,b_n)$ in the sense of Definition \ref{type}
below; this is why we use this ambiguous notation).
Thus if we take $R=\{b_k\mid k=1,\dots,n\}$,
for all $i\gg 0$ for given $\tau>0$, we
satisfy all assumptions of Proposition \ref{transform} for any $(X_i,B_i)$ and
$(X_i,B)$ with $X_i,B_i,B$, and $P_i$ instead of $X,B,B'$, and $P$
respectively, except for the $\mathbb{R}$-Cartier property of $K_{X_i}+B$ in (3), and the lc property in (4).
The $\mathbb{R}$-Cartier property  will hold if for example $X$ is $\mathbb{Q}$-factorial.
Moreover, if $a< 1$, then each $a_i<1$ for $i\gg 0$,
and there exists a crepant extremal divisorial extraction
$Y_i\to X_i$
of an exceptional prime b-divisor $E_i$ with centre $P_i$ on $X_i$ and $a_i=\mld(P_i,X_i,B_i)=a(E_i,X_i,B_i)$
[\ref{log models}, Theorem 3.1].
The extration $Y_i$ is $\mathbb{Q}$-factorial too [\ref{log models}, Theorem 3.1],
and $K_{Y_i}+B_{Y_i}+(1-a)E_i$ is $\mathbb{R}$-Cartier where
$B_{Y_i}$ denotes the birational transform of $B_i$ on $Y_i$.

Thus in the $\mathbb{Q}$-factorial case and
under (4), for all $i\gg 0$,
$(X_i,B)$ is $a$-lc, and $K_{Y_i}+B_{Y_i}+(1-a)E_i$
is seminegative$/X_i$.
Cf. the proof of Proposition \ref{aa} below, and Step~8 in
the proof of Proposition \ref{Weak finiteness} where
we either assume (4) or, we assume ACC for lc thresholds and derive (4)
from that assumption.

Finally, note that we can derive (4) from
ACC for lc thresholds in dimension $d$ when $X$ is $\mathbb{Q}$-factorial.
Indeed, if (4) does not hold, then  possibly after passing to a subsequence, we can construct a strictly increasing sequence of boundaries
$B_i'$ such that  $B_i<B_i'<B$, and
 $\ldis(X_i,B_i')=0$. Essentially, $B_i'$ is constructed by taking an appropriate lc threshold.
Moreover, the multiplicities of $B_i'$ will
satisfies the dcc with finitely many accumulation
points in $R$. This contradicts ACC for lc thresholds in dimension $d$.

\end{exa}

\begin{proof}(of Proposition \ref{aa})

Suppose that $K_{Y}+B_Y'+(1-a)E$ is not seminegative.
Then by property (1) of the proposition and the extremal property,
it is numerically positive$/X$.

On the other hand,
by (3) of Proposition \ref{transform},
$$
K_Y+B_Y'+\sum(1-a(E_i,X,B'))E_i\equiv 0/X,
$$
where by  Proposition \ref{transform} the discrepancy $a(E_i,X,B')\ge a$
for each prime component $E_i$ of $E=\sum E_i$.
Thus the $\mathbb{R}$-Cartier divisor
$$
\sum(a-a(E_i,X,B'))E_i=(K_Y+B_Y'+\sum(1-a(E_i,X,B'))E_i)-
(K_Y+B_Y'+(1-a)E)
$$
is numerically negative$/X$.
According to Negativity [\ref{log flips}, 1.1],
the divisor is effective and $\not=0$, that is,
each $a-a(E_i,X,B')>0$, a contradiction.
$\Box$
\end{proof}

The following result is the big chunk of (i) in
our Main Theorem \ref{result}, and it gives another application of Main Proposition when the support
of $B$ is not universally bounded.

\begin{prop}\label{mld-th} $\ACC$ for mld's and lc thresholds in dimension d implies $\ACC$ for $a$-lc
thresholds in the same dimension for all $a>0$, in particular, for canonical thresholds.
\end{prop}

\begin{proof}
Suppose that we have a monotonic increasing sequence
$t_i$ of $d$-dimensional $a$-lc thresholds, that is,
there exists a sequence $(X_i,B_i)$ of
$d$-dimensional log pairs with
boundaries $B_i\in \Gamma$, and
$\mathbb{R}$-divisors $M_i\in S$ on $X_i$ such that

(1) $(X_i,B_i)$ is $a$-lc;

and

(2) $t_i=\thr_a(M_i,X_i,B_i)$;

in particular,
on each $X_i$ there exists a point $P_i\in \Supp M_i\subset X_i$ of codimension $\ge 2$,
at which

(3) $\mld(P_i,X_i,B_i+t_iM_i)=a$; or

(4) $B_i+t_iM_i$ has a reduced component
(the payment for $\ldis$ in codimension $\ge 2$;
see see Remark \ref{th-rem}).

We need to verify the acc for the sequence $t_i$,
that is, the sequence stabilizes.

If for infinitely many $i$,  $B_i+t_i M_i$ has a reduced component $D_{i,j}$
as in (4), that is, $b_{i,j}+t_im_{i,j}=1$ for multiplicities in $D_{i,j}$, then $t_i=(1-b_{i,j})/m_{i,j}$ and stabilizes by the dcc for multiplicities
$b_{i,j}$ and $m_{i,j}$. This gives the acc in the case (4).

Thus after taking a subsequence, we can assume (3) for all $i$.
Note that by the lc property of $(X,B_i+t_iM_i)$ the limit
$t=\lim_{i\to \infty} t_i$ exists
because $t_i$ are bounded from above:
  $t_i\le \frac{1}{m_0}$ where $m_0=\min\{m\in S\}$; $t\ge 0$.
  We can apply Proposition \ref{transform} for each $X=X_i,B=B_i,B'=B_i+tM_i$, and $P=P_i$.
  Indeed, (1) of the proposition holds because $\dim X= d$.
The assumption (2) of the proposition follows from construction, in particular,
the multiplicities of $B'$ are as $b_{i,j}+tm_{i,j}$ and satisfy the dcc as their components $b_{i,j}$ and $m_{i,j}$ do.

The assumption (3) of the proposition
$\mld(P_i,X_i,B_i)\ge a$ holds by (1) above; $K_X+B'$ is $\mathbb R$-Cartier because each $M_i$ is $\mathbb R$-Cartier.

Finally, the assumption (4) of the proposition, that is, $(X,B')$ is lc
in a neighborhood of $P$, follows from $\ACC$ for lc thresholds.
Indeed, if the lc property does not hold for $i\gg 0$, then we get an increasing set of $t_i'=\lct(M_i,X_i,B_i)$ for infinitely many $i$, such that $t_i\le t_i'<t$ (The lc property in codimension $1$ holds
by construction.). This contradicts ACC for lc thresholds.

Therefore, by Proposition \ref{transform} and (1-3) $t=t_i$ for all $i\gg0$ and
$t_i$ stabilizes.

$\Box$
\end{proof}

Proposition \ref{transform} also gives some relations between different
$\ACC$ versions besides the ones for mld's and thresholds in the Introduction.
Now we recall some of them.

\begin{defn}[cf. {[\ref{K$^+$}, Section 18].}]\label{type}

(i) The {\em type order\/} is a direct sum of $\mathbb{R}$
countably many times, that is, the set of
sequences $(b_1,\dots,b_n)$ with $b_i\in \mathbb{R}, n\ge 0$, and
the following order:
$(b_1,\dots,b_m)\le (b_1',\dots,b_n')$ if either $n<m$ or
$n=m$ and each $b_i\le b_i'$.
The {\em maximal element} is the empty sequence with $n=0$.

(ii) A {\em type\/} of an $\mathbb{R}$-divisor $D=\sum d_i D_i$ on $X$,
where $D_i$ are distinct prime divisors on $X$,
is the sequence $(d_1,\dots,d_n)$ of its nonzero multiplicities
(in any possible ordering).
We usually do not think of $D$ with a specific
ordering of the prime components in mind, so $D$
can have several types.
Even one can add finitely many zeros.

(iii) (Cf. [\ref{K$^+$}, Definition~18.3].)
A log pair $(X,H+D)$ with $\mathbb{R}$-divisors $H,D$ and
prime divisors $D_1,\dots,D_n$ on $X$ has {\em maximal
$a$-lc type $(d_1,\dots,d_n)$ near\/} $Z\subset X$ and
respectively {\em at\/} $P\in X$ if $D=\sum d_i D_i$,
in particular, has
type $(d_1,\dots,d_n)$, if $(X,D)$ is $a$-lc near $Z$ and
respectively at $P$, and $(X,D')$ is not $a$-lc near $Z$ and
respectively at $P$ for any $\mathbb{R}$-divisor $D'=\sum d_i' D_i$
of type $(d_1',\dots,d_n')$
such that  $K+H+D'$ is $\mathbb{R}$-Cartier and
$D'>D$ in any neighborhood of $Z$ and respectively of $P$.
For $H=0$, the pair $(X,D)$ has that maximal property.

(iv) (Cf. [\ref{K$^+$}, 18.15.1].) $\mathcal S_d\text{(Fano)}$ is the set
of types $(b_1,\dots,b_n)$ such that there is
a nonsingular Fano variety $X$ of dimension at most $d$ and
a boundary $B$ of type $(b_1,\dots,b_n)$ such that $K+B\equiv 0$.
(See Example \ref{acc-exa}, (1) below.)

(v) (Cf. [\ref{K$^+$}, 18.15.$\overline{1}$].) $\mathcal S_d\text{(global)}$ is
the set of types $(b_1,\dots,b_n)$ such that there is
a proper normal variety $X$ of dimension at most $d$ and
a boundary $B$ of type $(b_1,\dots,b_n)$ such that $(X,B)$ is lc,
and $K+B\equiv 0$.
(See Example \ref{acc-exa}, (2) below.)
We denote by $\mathcal S_d^0\text{(global)}$ its subset with
nonklt $(X,B)$.

(vi) (Cf. [\ref{K$^+$}, 18.15.$\overline{2}$].)
$\mathcal S_{a,d}\text{(local)}$ is the set of
types $(b_1,\dots,b_n)$ such that there is a pointed
$\mathbb{Q}$-factorial variety $P\in X$ of dimension at most $d$, and
prime divisors $D_1,\dots,D_n$ on $X$ such that
$B=\sum b_i D_i$ is a boundary,
$$
P\in\cap D_i \text{ and }
\mld(P,X,B)\ge a,
$$
actually $=a$ or $n\ge 1$ and of maximal lc type in some point of any neighborhood of $P$ (see Example \ref{acc-exa}, (4) below),
or equivalently,
$B=\sum b_i D_i$ is a boundary, and
$(X,B)$ locally has maximal $a$-lc
type $(b_1,\dots,b_n)$ at $P$ with given divisors $D_i$.
(See Example \ref{acc-exa}, (3) below.)

(vii) (Cf. [\ref{K$^+$}, 18.15.$\overline{2}$].)
$\overline{\mathcal S}_{a,d}\text{(local)}$ is the set of
types $(b_1,\dots,b_n)$ such that there is
a $\mathbb{Q}$-factorial variety of dimension at most $d$,
prime divisors $D_1,\dots,D_n$ on $X$,
and a closed subvariety $Z\subset X$ such that
$B=\sum b_i D_i$ is a boundary,
$(X,B)$ is $a$-lc near $Z$,
every $D_i$ intersects $Z$, and $(X,B)$ has maximal $a$-lc
type $(b_{i_1},\dots,b_{i_l}), 1\le i_1<\dots<i_l\le n$,
at some of the intersection points of $D_i$  with the divisors $D_{i_1},\dots,D_{i_l}$
passing through such a point, or
equivalently,
$B=\sum b_i D_i$ is a boundary, and
$(X,B)$ locally has maximal $a$-lc type $(b_1,\dots,b_n)$ near $Z$
with the divisors $D_i$.

(viii) (Cf. [\ref{K$^+$}, 18.15.3].)
$\mathcal S_{a,d}^0\text{(local)}$ is the set of
types $(b_1,\dots,b_n)$ such that there is a pointed
$\mathbb{Q}$-factorial variety $P\in X$ of dimension at most $d$, and
prime divisors $D_1,\dots,D_n$ on $X$ such that
$B=\sum b_i D_i$ has a reduced component,
$$
P\in\cap D_i \text{ and }
\mld(P,X,B)\ge a,
$$
actually $=a$ or $n\ge 1$ and of maximal lc type in some point in any neighborhood of $P$ (see Example \ref{acc-exa}, (4) below),
or equivalently,
$B=\sum b_i D_i$ has a reduced component, and
$(X,B)$ locally has maximal $a$-lc
type $(b_1,\dots,b_n)$ at $P$ with given divisors $D_i$.

(ix) (Cf. [\ref{K$^+$}, 18.15.$\overline{3}$].)
$\overline{\mathcal S}_{a,d}^0\text{(local)}$ is the set of
types $(b_1,\dots,b_n)$ such that there is
a $\mathbb{Q}$-factorial variety $X$ of dimension at most $d$,
a subset $Z\subset X$, and
prime divisors $D_1,\dots,D_n$ on $X$ such that
$B=\sum b_i D_i$ has a reduced component,
$(X,B)$ is $a$-lc near $Z$,
$Z$ is in the reduced part of $B$,
every $D_i$ intersects $Z$, and $(X,B)$ has maximal $a$-lc
type $(b_{i_1},\dots,b_{i_l}), 1\le i_1<\dots<i_l\le n$,
at some of the intersection points of $D_i$ with the divisors $D_{i_1},\dots,D_{i_l}$
passing through such a point, or
equivalently,
$B$ has a reduced component,
$Z$ is in the reduced part of $B$,
$(X,B)$ locally has maximal $a$-lc type $(b_1,\dots,b_n)$ near $Z$
with given divisors $D_i$.

(x) (Cf. [\ref{K$^+$}, 18.15.1].) $\mathcal S_d\text{(Mori-Fano)}$ is the set of
types $(b_1,\dots,b_n)$ such that there is
a Mori-Fano variety $X$ of dimension at most $d$, and
a boundary $B$ of type $(b_1,\dots,b_n)$ such that $(X,B)$ is lc,
and $K+B\equiv 0$.
We denote by $\mathcal S_d^0\text{(Mori-Fano)}$ its subset with
nonklt $(X,B)$.

(xi) (Cf. [\ref{K$^+$}, 18.15.1].) $\mathcal S_d\text{(Mori-Fano cn)}$
is the set of
types $(b_1,\dots,b_n)$ such that there is
a cn Mori-Fano variety $X$ of dimension at most $d$ and
a boundary $B$ of type $(b_1,\dots,b_n)$,
and $K+B\equiv 0$.

We consider each of the above sets $\mathcal S_d\text{(Fano)},\dots,
\mathcal S_d\text{(Mori-Fano cn)}$ as a subset of the order $\mathcal B$.
Thus it has ordering induced from $\mathcal B$.

For $a=0$, we set $\mathcal S_d=\mathcal S_{a,d}$, e.g.,
$\overline{\mathcal S}_d^0\text{(local)}=
\overline{\mathcal S}_{0,d}^0\text{(local)}$.
Some of them are slightly more general
than in [\ref{K$^+$}, 18.15].
 \end{defn}

Nonetheless we expect the same.

\begin{conj}[cf. {[\ref{K$^+$}, Conjecture~18.16]}]\label{acc-conj}
Each set $\mathcal S_d\text{(global)}$, $\mathcal S_d^0\text{(global)}$,
$\mathcal S_{a,d}\text{(local)}$, $\overline{\mathcal S}_{a,d}\text{(local)}$, $\mathcal S_{a,d}^0\text{(local)}$,
$\overline{\mathcal S}_{a,d}^0\text{(local)}$, $\mathcal S_d\text{(Mori-Fano)}$, $\mathcal S_d\text{(Mori-Fano cn)}$ satisfies the acc.
\end{conj}

\begin{rem}\label{acc-rem}
(1)
Equivalently, in Definition \ref{type} (iii) for effective $D$
and $D'>D$,
the assumption $D'=\sum d_i' D_i$ has type $(d_1',\dots,d_n')$
can be replaced by the strict inequality of types:
$D'$ has type $(d_1',\dots,d_m')>(d_1,\dots,d_n)$
when each $d_i\not=0$, that is, $m=n$ in this situation.

Of course, the maximal $a$-lc property depends on $a,H$ and
divisors $D_i$.
Actually, at $P$, it depends only on $D$ itself;
for $n=0$, $D=0$ and $D'>0$ with $\Supp D'=0$ does not
exist (cf. Example \ref{acc-exa}, (4) below).
However if we consider a type $(d_1,\dots,d_n)$ of $D$ with
all $d_i\not=0$, maximal $a$-lc near $Z$ also depends only
on $D$ itself.
This condition can be stated as the maximal type
$(d_1,\dots,d_n)$ for $D$.
It is unique up to permutation.

In general, we can add some $b_i=0$ (see Example \ref{acc-exa} for $n=2$, or proof of ??\footnote{?} below)

(2) In Definition \ref{type}, (vi-ix) the $\mathbb{Q}$-factorial assumption
can be replaced by the $\mathbb{Q}$-Cartier property of prime
divisors $D_i$ (cf. Example \ref{acc-exa}, (3) below).

(3) In Definition \ref{type}, (vii) and (ix) for $a>0$,
$(X,B)$ is plt near $Z$.
Thus the reduced part $B_0=\sum_{b_i=1} D_i$ of $B$
is locally irreducible near each point of $Z$.
Actually, in this situation $a\le 1$
in dimension $d\ge 2$.

(4)
Since we omit the lc assumption in Definition \ref{type}, (iv) and (xi),
the sets
$S_d\text{(Fano)}$ and $S_d\text{(Mori-Fano cn)}$ are not
subsets of $S_d\text{(Mori-Fano)}$, and
their intersections with $S_d\text{(Mori-Fano)}$ are
determined by the lc property of $(X,B)$.
Nonetheless acc is known for $S_d\text{(Fano)}$ (see
Example \ref{acc-exa}, (1) below) and
expected for $S_d\text{(Mori-Fano cn)}$.
ACC for $S_d\text{(Mori-Fano cn)}$ follows
from Conjecture \ref{weak-BAB} in dimension $d$
as so does acc for $S_d\text{(Fano)}$ from the boundedness
of nonsingular Fano varieties in Example \ref{acc-exa}, (1).
Moreover, we can omit in both cases the assumption $b_i\le 1$, and
in the last case the assumption
$\rho(X)=1$ because the boundedness of cn Fano varieties
in any dimension $d$ is expected (cf. Example \ref{acc-exa}, (5) below).
\end{rem}

\begin{exa}\label{acc-exa}
(1)
$\mathcal S_d\text{(Fano)}$ satisfies the acc in
any dimension $d$ by the boundedness of nonsingular
Fano varieties $X$ of dimension $\le d$ [\ref{KM}].
Indeed, there exists a generic curve $C\subset X$ which
positively intersects each prime divisor $D_i$ on $X$ and
with bounded $(-K\cdot C)$: $C$ is
a generic curve section for an embedding $X\subset \mathbb{P}^N$ of bounded degree in a fixed
projective space $\mathbb{P}^N$.
Then for positive integers $m_i=(D_i\cdot C)$,
$\sum b_i m_i=(B\cdot C)=(-K\cdot C)$.
On the other hand,
for any increasing sequence of types $(b_1^l,\dots,b_{n_l}^l),l=1,2, \dots$,
we can suppose that their sizes stabilize: $n_l=n$ for all $l\gg 0$, and
each $b_i^l>0$.
Therefore, after taking a subsequence,
the multiplicities $m_i^l=(D_i^l\cdot C)$ stabilize too:
for each $i=1,\dots,n$, $m_i^l=m_i>0$ for all $l\gg 0$.
Hence the types stabilize:
for each $i=1,\dots,n$, $b_i^l=b_i>0$ for all $l\gg 0$
(cf. the proof of [\ref{log flips}, Second Termination~4.9]).

Equivalently, for a finite set $R$ of real (nonnegative) numbers,
we can consider the set of types $(d_1,\dots,d_n)$ such that
each $d_i\ge 0$ and $\sum d_i=d\in R$
(cf. Example \ref{d=1} below).
For example, the case $R=\{2\}$ includes $\mathcal S_1\text{(Fano)}$
with an extra condition, namely, each $d_i=b_i\le 1$.

According to our arguments for the acc,
the assumption that $(X,B)$ is lc,
and other ones on singularities of the log pair
are not necessary, in particular, we can omit
the assumption $b_i\le 1$.
Hence acc holds for $\mathcal S_d\text{(fano)}$ as in
[\ref{K$^+$}, 18.15.1].

Let $X=\mathbb{P}^d$ be the projective space of dimension $d$, and
$D$ a generic hypersurface in $\mathbb{P}^d$ of degree $d+2$.
Then
$$
K_{\mathbb{P}^d}+\frac{d+1}{d+2} D\equiv 0.
$$
Thus $((d+1)/(d+2))\in \mathcal S_d\text{(Fano)}$, and
the dimension condition for sets in Definition \ref{type}, (iv-v), and (x-xi)
is necessary to satisfy the acc in Conjecture \ref{acc-conj}.
Similarly, for all other sets in the conjecture.

(2) However, for $\mathcal S_d\text{(global)}$ and
$\mathcal S_d\text{(Mori-Fano)}$,
the assumption that $(X,B)$ is lc
is very important.
Let $Q_n\subset \mathbb{P}^{(n+1)}$ be the cone over
a rational normal curve of degree $n$ with
a line generator $L$.
Then for a generic hyperplane section $H$,
$$
-K=(n+2)L\equiv 3L+\frac {n-1}n H.
$$
If we replace $3L$ by $L_1+L_2+L_3$ or
$L_1/2+\dots +L_6/2$ with distinct generators $L_i$,
we construct strictly increasing sequences of types
$(1,1,1,(n-1)/n)$ and $(1/2,1/2,1/2,1/2,1/2,1/2,(n-1)/n)$
respectively.
However, they are not in $\mathcal S_d\text{(global)}$
because $(Q_n,L_1+L_2+L_3+(n-1)H/n)$ and
$(Q_n,L_1/2+\dots+L_6/2+(n-1)H/n)$ are not lc
(at the vertex of $Q_n$).

(3) The $\mathbb{Q}$-factorial property in Definition \ref{type},(vi-ix)
is very important, too
 (cf. Remark \ref{acc-rem}, (2) above).
Let $f\colon Y\to X$ be a contraction of
a nonsingular rational curve $C$ on a nonsingular $3$-fold
$Y$, and $D_1,D_2$ two nonsingular prime divisors on $Y$
with intersection only along $C$ with normal
crossings.
Set $-n=C^2$ on $D_1$.
For any $n\ge 2$ there exists such a contraction, e.g., toric one.
Then $K+B\equiv 0/X$ for $B=D_1+(n-2)D_2/n$.
Thus we have a strictly increasing sequence of types
$(1,(n-2)/n)$ which does not belong (entirely) to any set in Conjecture \ref{acc-conj}
if it satisfies the acc.
That is in Definition \ref{type}, (v) the {\em proper\/} assumption is
very important.
The same types correspond to the image $(X,f_*B)$.
However it does not belong to the sets in Definition \ref{type}, (vi-ix)
because $X$ is not $\mathbb{Q}$-factorial (cf. Remark \ref{acc-rem}, (2) above).

Let $D_3$ be a divisor which transversally
intersects $D_1,D_2$ in a single point.
(Again such a divisor exists in a toric case.)
Then, for
$B'=D_1+D_3+(n-1)/nD_2$, $K+B'\equiv 0/ X$, and
$(X,f_*B')$ is exactly lc near $P=f(C): \ldis(X,f_*B')=0$,
but $(n-1)/n$ does not satisfies the acc and
is not a counter example to $\ACC$ for thresholds
since $f_*D_2$ is not a $\mathbb{Q}$-Cartier divisor.
Similar examples can be constructed for
any $a$ instead of $0$.

However it is expected that (the existence of
$\mathbb Q$-factorialization implies that),
for any strictly increasing types, $\ldis$ $0$
is never attained at $P$.

(4) If in Definition \ref{type}, (vi) and (viii) $n=0$, that is,
the type itself is {\em maximal\/},
then $B=0$, and $\mld(P,X,0)\ge a$ only but can be $\not=a$.
More generally, $(X,H+0)$ has maximal $a$-lc type $\emptyset$
near $Z$ and respectively at $P$ if locally
$K+H$ is $\mathbb{R}$-Cartier and $(X,H)$ is $a$-lc near $Z$ and
respectively $\mld(P,X,H)\ge a$ but not necessary $=a$.

(5)
Acc of $\mathcal S_d\text{(Mori-Fano cn)}$ holds
for $d\le 3$ by [\ref{KMMT}].
Moreover, one of the boundary properties of $B$, namely
that each $b_i\le 1$, and
the condition $\rho(X)=1$ are not necessary (cf. Remark \ref{acc-rem}, (4) above).
\end{exa}

The main result of this section is

\begin{thm}\label{main-th}
$\ACC$ for mld's and lc thresholds in dimension $\le d$ imply

(i)
Acc for $\mathcal S_{a,d}\text{(local)},
\overline{\mathcal S}_{a,d}\text{(local)},
\mathcal S_{a,d}^0\text{(local)}$, and
$\overline{\mathcal S}_{a,d}^0\text{(local)}$ with
any $a>0$;

Acc for $\mathcal S_d\text{(Mori-Fano cn)}$, $\LMMP$ and
$\ACC$ for mld's in dimension $\le d$ imply

(ii) acc for $\mathcal S_{d}\text{(global)}$, and
$\mathcal S_{d}\text{(Mori-Fano)}$;

(iii) acc for $\mathcal S_{d+1}^0\text{(Mori-Fano)}$;

(iv)
acc for $\mathcal S_{d+1}\text{(local)},
\overline{\mathcal S}_{d+1}\text{(local)},
\mathcal S_{d+1}^0\text{(local)}$, and
$\overline{\mathcal S}_{d+1}^0\text{(local)}$; and

(v) $\ACC$ for
lc thresholds in dimension $\le d+1$;

If in addition, $\LMMP$ holds in dimension $d+1$,
then

(vi) acc for $\mathcal S_{d+1}^0\text{(global)}$ holds.

\end{thm}

\begin{adden}\label{main-th-adden} Acc for $\mathcal S_d\text{(Mori-Fano cn)}$ can be replaced by
Conjecture \ref{weak-BAB} in dimension $\le d$ (everywhere!)
because the latter implies
the former.
(Cf. Example \ref{acc-exa}, (1) and Remark \ref{acc-rem}, (4) above.)
\end{adden}

\begin{thm}\label{main-th2}
Let $\Gamma\subset [0,1]$ be a set satisfying the dcc.
Then Theorem \ref{main-th} holds when each $B\in \Gamma$.
\end{thm}

Notation: We denote the corresponding subsets by
$\mathcal S_d\text{($\Gamma$, local)}$, etc. Of course, we can apply it in other cases too:
e.g., $\mathcal S_d\text{($\Gamma$, Fano)}$,
and similar results hold in these cases.

\begin{adden}\label{add-main-th2}
 Moreover, then $\Gamma$ can be assumed to be finite, that is,
there exists its finite subset $\Gamma_f$ such that
$
\mathcal S_{a,d}\text{($\Gamma$, local)}=\mathcal S_{a,d}\text{($\Gamma_f$, local)}$, $ \mathcal S_{d}\text{($\Gamma$, global)}=\mathcal S_{d}\text{($\Gamma_f$, global)}$, etc.

\end{adden}

We also have  $\mathcal S_d\text{($\Gamma$, Fano)}=\mathcal S_d\text{($\Gamma_f$, Fano)}$.

\begin{cor}[ACC for anticanonical (ac) thresholds]\label{anticanonical}
Assume acc for $\mathcal S_d\text{(Mori-Fano cn)}$, $\LMMP$ and
$\ACC$ for mld's in dimension $\le d$. Then,

(i) acc for $\mathcal S_{a,d}\text{(global FT)}$ where the $\mathcal S_{a,d}\text{(global FT)}$ is the set of types $(b_1,\dots,b_n)$
such that there is  a FT variety $X$ of dimension at most $d$, an ample Cartier
divisor $H$ on $X$, and
a boundary $B$ of type $(b_1,\dots,b_n)$ such that $(X,B)$ is lc,
and $K+B+aH\equiv 0$; in particular,
$\mathcal S_{0,d}\text{(global FT)}=\mathcal S_d\text{(global FT)}$ the subset in $\mathcal S_d\text{(global)}$ corresponding to FT $X$;

(ii)
for subsets $\{0\}\subseteq\Gamma\subseteq [0,1]$ and
$\{0\}\subsetneq S\subseteq \mathbb{R}$ of nonnegative numbers
satisfying the dcc,\\

{\rm (ACC)}
the following subset of $\mathbb{R}^+\cup\{+\infty\}=\{r\in \mathbb{R}\mid r\ge 0\}\cup \{+\infty\}$
$$\overline{\mathcal{A}}_{a,d}(\Gamma,S)=\{\act_a(M,X,B)|\mbox{$X$ is complete, $\dim X=d$, $B\in \Gamma$, $(X,B)$ is lc},$$
$$\mbox{$X$ is FT, $K+B$ is seminegative and $M$ is an S-Cartier divisor on $X$}\}.$$

satisfies  the acc; $+\infty$ corresponds to the case $M=0$,
where $t=\act_a(M,X,B)$  means that $K+B+tM+aH\equiv 0$ on $X$
for some ample Cartier divisor $H$ on $X$, and the $S$-Cartier property means
that $M$ is a linear combination of ample Cartier divisors with
multiplicities in $S$.

In particular, the $\ACC$ holds for $a=0$ and $S=\{1\}$  which gives
the anticanonical threshold (see [\ref{Isk-Sh}, p.~47]).

(iii) the log Fano indices, that is, a maximal
real positive number $a$ such that
$K+B+aH\equiv 0$ for ample Cartier
divisors $H$, satisfies the acc for
the lc pairs $(X,B)$, with FT variety $X$ of dimension $\le d$
with $B\in \Gamma$ as in (ii).

\end{cor}

\begin{adden}\label{act-adden} Acc for $\mathcal S_d\text{(Mori-Fano cn)}$ can be replaced by
Conjecture \ref{weak-BAB} in dimension $\le d$.
\end{adden}

\begin{rem}
We expect that Corollary \ref{anticanonical} holds when FT is omitted, that is,
for $a>0$, $(X,B)$ is just a lc Fano variety as in Theorem \ref{main-th} (ii).

\end{rem}

\begin{cor}\label{act-cor} Let $\Gamma\subset [0,1]$ be a set satisfying the dcc.
Then Corollary \ref{anticanonical} (i) holds when each $B\in \Gamma$.
\end{cor}

\begin{adden}\label{act-adden2} Moreover, then $\Gamma$ can be assumed to be finite in Corollary \ref{act-cor} (i), that is,
there exists its finite subset $\Gamma_f$ such that
$\mathcal S_{a,d}\text{($\Gamma$, global)}=
\mathcal S_{a,d}\text{($\Gamma_f$, global)}$.
\end{adden}

\begin{rem} If $\rho (X)=1$, then $a$-{\em anticanonical\/} ({\em $a$-ac\/})
{\em threshold\/} is well defined for any ample Cartier divisor $H$, e.g., for
such a generator in $\Pic(X)$:
there exists  a (unique) real number $t$ such that
$$
K+B+tM+aH\equiv 0.
$$
If $K+B+aH$ is seminegative then $t\ge 0$. The condition $\rho(X)=1$ replaces the
$\mathbb{Q}$-factorial property of the local case. For $a=0$, we get the ac threshold [\ref{Isk-Sh}, p.~47].

In Corollaries \ref{anticanonical} and \ref{act-cor} we can suppose that $a$
is varying in a dcc set. Then it is expected that the corresponding
thresholds $t$ in dimension $\le d$ satisfies
the acc (ACC conjecture). This is clear from the proof of Corollaries \ref{anticanonical} and \ref{act-cor} (below).
\end{rem}

\begin{proof} (of Corollary \ref{anticanonical})

The case (i) follows from its counterpart in Theorem \ref{main-th} (ii) for $a=0$. To apply the theorem we replace the boundary $B$ with  $B+aH$ with an appropriate choice of $H$ (see proof of (ii) below). The type of $B$ will be extended by that of $aH$. Since the latter has finitely many possible multiplicities, acc for $B$ is equivalent to acc of
extended types to which we apply Theorem~1.10.

(ii) Suppose that such thresholds do not satisfy the acc. Let $\Omega$ be an infinite set of such thresholds, which satisfies the dcc. Now take a $t\in \Omega$. Also take $X,B,H$ and $M\neq 0$ corresponding to $t$. Since $M$ is $S$-Cartier, there are $s_j\in S$ and ample Cartier divisors $H_j$ such that $M=\sum_js_jH_j$. By anticanonical boundedness, $ts_j$ is bounded. By effective base point freeness [\ref{ebf}], there is $h$, a natural number (not depending on $X, H,H_j$ but depending  only on the
dimension $d$),
such that $hH_j$ and $hH$ are free divisors and $(a/h),(ts_j/h)\in [0,1]$.

 Write $M\equiv \sum_j(ts_j/h)H_j'$ where $H_j'\in |hH_j|$ is a general member. For $d\ge 2$, general $H_j'$ is irreducible. For $d=1$, the number of components in $H_j'$ is bounded. Thus we can assume that
$$
(X,B+  \sum_j(ts_j/h)H_j'+(a/h)H')
$$
is lc where  $H'\in |hH|$ is general. In particular, the possible multiplicities of $B+  \sum_j(ts_j/h)H_j'+(a/h)H'$ satisfies the dcc. Now use Addendum \ref{add-main-th2} for $S_d(\Gamma_f,global)$.

(iii) If $S={1}$ then any $M$ is an ample Cartier divisor $H$ and any ac threshold
satisfies $K+B+tH\equiv 0$. Thus possible $t$ satisfy the acc by (ii).
This implies acc for the Fano indices.

$\Box$
\end{proof}

\begin{proof}(of Addnedum \ref{act-adden})

Same as Addendum \ref{main-th-adden} below. $\Box$
\end{proof}

\begin{proof}(of Corollary \ref{act-cor})

There is a natural number $h$ such that $hH$ is free where $h$ does not depend on $X,H$. Choose $h$ big enough such that $a/h\in [0,1]$. Now for a general $H'\in |hH|$, $K+B+(a/h)H'\equiv 0$ and $(X,B+(a/h)H')$ is klt. Since $B\in \Gamma$, possible multiplicities of $B+(a/h)H'$ satisfy the dcc because $\Gamma$ satisfies the dcc, $a,h$ are fixed and $H'$ is a reduced Cartier divisor. Therefore, the result follows from acc for $\mathcal S_{d}\text{(global)}$ as in Theorem \ref{main-th2} (ii).

$\Box$
\end{proof}

\begin{proof}(of Addendum \ref{act-adden2})

We can use the extension of boundaries $B$ as
in the proof of Corollary \ref{anticanonical}, and then use
Addendum \ref{add-main-th2}.
$\Box$
\end{proof}

\begin{proof}(of Theorem \ref{main-th})
Each statement follows from the same statement under
the assumption $B\in \Gamma$ for some $\Gamma$ under
the dcc.
Thus it follows from the corresponding statement in Theorem \ref{main-th2}.

For the statements (iii) and (vii) such a set $\Gamma$
is given by assumptions.

In the other cases, we need to verify that each increasing sequence of
types $(b_1^l,\dots,b_{n_l}^l), l=1,2,\dots,$
stabilizes.
By definition of the ordering, their sizes stabilize:
$n_l=n$ for all $l\gg 0$.
Then for the corresponding pairs $(X_l,B_l)$,
$B_l=\sum_{i=1}^{n_l} b_i^l D_i^l\in \Gamma$
where
$$
\Gamma=\{b_i^l\mid i=1,\dots,n,\text{ and }
l=1,2,\dots\}=\cup_{i=1}^n \Gamma_i, \text{ and }
\Gamma_i=\{b_i^l\mid l=1,2,\dots\}.
$$
Since each sequence $b_i^l,l=1,2,\dots,$ increases,
the sets $\Gamma_i$ and $\Gamma$ satisfy the dcc.
\end{proof}

\begin{proof}(of Addendum \ref{main-th-adden}) Immediate by Remark \ref{acc-rem}~(4).
$\Box$
\end{proof}

\begin{proof}(of Theorem \ref{main-th2})  (i) By inclusions
$$
\mathcal S_{a,d}^0\text{($\Gamma$, local)}\subset
\mathcal S_{a,d}\text{($\Gamma$, local)} \text{ and }
\overline{\mathcal S}_{a,d}^0\text{($\Gamma$, local)}\subset
\overline{\mathcal S}_{a,d}\text{($\Gamma$, local)}
$$
It is enough to prove acc for the ambient sets.
On the other hand,
it is known (almost by definition) that,
for any type $(b_1,\dots,b_n)$ in
$\overline{\mathcal S}_{a,d}\text{($\Gamma$, local)}$
and any of its component $b_i$ there exists a type
$(b_1',\dots,b_m')$ in $\mathcal S_{a,d}\text{($\Gamma$, local)}$
with a component $b_{i'}'=b_i$: take a type corresponding to
an appropriate intersection point of $D_i\cap Z$ (cf. [\ref{K$^+$}, 18.19.1]).
Thus it is enough to prove acc for $\mathcal S_{a,d}\text{($\Gamma$,local)}$.

Let $(b_1^l,\dots,b_{n_l}^l),l=1,2,\dots,$ be
an increasing sequence of types in the set $\mathcal S_{a,d}\text{($\Gamma$, local)}$. The fact that this sequence of types is increasing, implies that $n_l$ is bounded and it stabilizes: $n_l=n$ for $l\gg 0$. So we can assume that $n_l=n$ for any $l$. Since the sequence stabilizes for $n=0$, we
can suppose that $n\ge 1$.
According to construction we have a sequence of
pointed $\mathbb{Q}$-factorial varieties $P_l\in X_l$ of dimension
 $d$ and prime divisors $D_1^l,\dots,D_n^l$ on $X_l$
such that $B_l=\sum b_i^l D_i^l$ is a boundary,
$$
P_l\in\cap D_i^l \text{ and }
\mld(P_l,X_l,B_l)=a.
$$
and $P_l$ has codimension $\ge 2$, or $(b_1^l,\dots,b_{n_l}^l)$ is maximal lc at some point $Q_l$ near $P_l$
(possibly of codimension $1$ but not a closed point since $a>0$). Note that in the last case, all $D_i^l$ with $b_i^l<1$ pass through
$Q_l$ by the maximal lc property. Thus a subvector of $(b_1^l,\dots,b_{n_1}^l)$ and so the set of all $(b_1^l,\dots,b_{n_1}^l)$
satisfies the acc by ACC for lc thresholds in dimension $\le d-1$ (see the arguments below). Thus taking a subsequence we can suppose that the
first case: $\mld(P_l,X_l,B_l)=a$.

We can choose  a
subsequence such that the limits below exist (e.g., unique) by monotonic increasing and
boundedness $(\le 1$, see Example \ref{bb} above)
$$
b_i=\lim_{i\to\infty}b_i^l ~~\mbox{for}~~i=1,\dots,n,~~~ B'=\sum_{i=1}^n
b_iD_i^l,~~\text{ and }
~~~R=\{b_i\mid i=1,\dots,n\}
.
$$

Then for any $\tau>0$,
$\parallel B_l-B'\parallel<\tau$ for all $l\gg 0$.
Note that $K_{X_l}+B'$ is $\mathbb{R}$-Cartier because
$X_l$ is $\mathbb{Q}$-factorial.
By ACC for lc thresholds and Proposition \ref{transform}
for $X=X_l,B=B_l,P=P_l,$ and every $l\gg 0$,
we can assume that $(X_l,B')$ is lc near $P_l$, and
$a$-lc at $P_l$; $a>0$ by assumptions. Therefore,
$\mld(P_l,X_l,B')\ge a=\mld(P_l,X_l,B_l)$.

We can derive the lc property, (4) of Main Proposition \ref{transform},
of $(X_l,B')$ from the assumptions as follows (cf. proof of Proposition \ref{mld-th}). If $(X_l,B')$ is not lc near $P_l$  for $l\gg 0$, then (since $X_l$ is $\mathbb{Q}$-factorial), for infinitely many $l$ there is $G_l=\sum_{i=1}^n g_{i}^lD_i^l $ such that  $B_l\le G_l\le B'$ and such that $(X_l,G_l)$ is precisely lc (i.e., lc but not klt) near $P_l$. The set of multiplicities of those $G_l$ satisfies the dcc and is not finite. We can assume that $\{g_{1}^l\}$ is not finite
but increasing and that $D_1^l$ contains a lc centre. So, $g_{1}^l$ is the lc threshold of $D_1^l$ with respect to $(X_l,G_l-g_{1}^lD_1^l)$. This contradicts ACC for lc thresholds.

Now by Monotonicity of mld's (see [Sh. 3-fold log flips,
1.3.3]) and since $B'\ge B_l$,
the sequence stabilizes: $B'=B_l$ for every $l\gg 0$.
This proves the acc.

(ii) This will be established in the weak finiteness section
modulo (v) in dimension $d$ which can be assumed by induction.

(iii) Let $(X_l,B_l)$ be a pair of dimension $d+1$ for each $l$ such that $B_l=\sum_{i=1}^{n_l} b_{i}^lD_{i}^l$ has a type
$(b_1^l,\dots,b_{n_l}^l) \in \mathcal S_{d+1}^0\text{($\Gamma$, Mori-Fano)}$ such that these types are strictly increasing with
respect to $l$. We can assume that $\{b_{1}^l\}$ is a strictly increasing sequence. By assumptions, $(X_l,B_l)$ is lc but not klt.  We can take a
strictly lt model  $(Y_l,B_{Y_l})$ for $(X_l,B_l)$; this needs special termination
and existence of flips in dimesion $d+1$ which follow from LMMP in dimension $d$ [\ref{F}].

Suppose that $D_{1}^l$ intersects $\LCS(X_l,B_l)$ for infinitely many $l$. Thus, the birational transform of $D_{1}^l$ intersects the
reduced part of $B_{Y_l}$. Then using adjunction, restrict to an appropriate component in the reduced part of the boundary which intersects
$D_{1}^l$. The multiplicities that we get are of the
following type
$$
1>b'=\frac{m-1}m+\sum \frac{c_i}mb_{i}^l
$$
with natural numbers $m,c_i$ [\ref{log flips}, 3.10].

\begin{lem}\label{restriction}
Any set of such $b'$ satisfies the dcc where $b_{i}^l\in \Gamma$
(cf. [\ref{log flips}, Second termination 4.9] and [\ref{K$^+$}, 18.21.4]).
Moreover, if it is finite, then the set of $b_{i}^l$ is finite.
\end{lem}

Here the finiteness of the set of $b'$ comes from induction on $d$ (part (ii)). So we get a contradiction.

Then, we can assume that $D_{1}^l$ does not intersect $\LCS(X_l, B_l)$. There is an extremal ray $R_l$ on $Y_l$ such that the birational
transform of $D_{1}^l$ intersects $R_l$ positively. If $R_l$ is of fibre type, then by restricting to the general fibre and using induction
on $d$ we get a contradiction. So assume otherwise. The reduced part of $B_{Y_l}$ intersects $R_l$, otherwise $R_l$ corresponds to a
flipping or divisorial type extremal ray $R_l'$ on $X_l$. This is not possible since $\rho(X_l)=1$.

Let $(Y^+_l,B^+_{Y_l})$ be the model after operating on $R_l$ (i.e., after a flip or divisorial contraction). Thus, the birational
transform of $D_{1}^l$ on $Y^+_l$ intersects the reduced part of $B^+_{Y_l}$. We get a contradiction as above by restricting to a
component of the reduced part of $B^+_{Y_l}$.

(iv)
As above it is enough to verify acc for $\mathcal
S_{d+1}\text{($\Gamma$, local)}$.
Suppose that there are $(X_l,B_l)$ such that  $B_l=\sum_{i=1}^{n_l} b_{i}^lD_{i}^l$ has a type  $(b_1^l,\dots,b_{n_l}^l)$ in
$\mathcal S_{d+1}\text{($\Gamma$, local)}$ such that these types are strictly increasing with respect to $l$. We can assume
that the set $\{b_{l}^l\}$
is strictly increasing. If for infinitely many $l$, $D_{1}^l$ passes through a lc centre of $(X_l,B_l)$
of dimension $\geq 1$, then by taking hyperplane sections, we reduce the problem to dimension
$\leq d$ for which we may assume that the theorem is already proved.

So, we assume that none of $D_{1}^l$ passes through a lc centre of $(X_l,B_l)$ of dimension
$\geq 1$. Now, take a strictly lt model of each $(X_l,B_l)$. Then using
adjunction, restrict to an appropriate exceptional divisor in the reduced part of the boundary
which intersects the birational transform of $D_{1}^l$. The multiplicities that we get are as in Lemma \ref{restriction}.
We get a contradiction by (ii).

(v) This is proved exactly as in (iv) using induction on $d$ (part (ii)).

(vi) Let $(X_l,B_l)$ be a $d+1$-dimensional pair for each $l$ such that $B_l=\sum_{i=1}^{n_l} b_{i}^lD_{i}^l$ has a type
$(b_1^l,\dots,b_{n_l}^l)\in \mathcal S_{d+1}^0\text{($\Gamma$,global)}$ such that these types are strictly increasing with respect
to $l$. We can assume that $\{b_{1}^l\}$ is a strictly increasing sequence. By assumptions, $(X_l,B_l)$ is lc but not klt.
As in the proof of Proposition \ref{Weak finiteness}, run the anti-LMMP on $D_{1}^l$. After finitely many steps, either we get
a fibration or the Mori-Fano case. For the former case we use
induction on $d$ and for the latter case use (iii).

 $\Box$
\end{proof}

To prove Addendum \ref{add-main-th2} we use the following

\begin{lem}\label{cc} Any suborder $\mathcal S\subset\mathcal B$
satisfies the acc if each $(b_1,\dots,b_n)$ in $\mathcal S$ is $\in R$,
that is, each $b_i\in R$,
for some fixed finite set of real numbers $R$.
The converse holds, that is, there exists
finite $R$ such that each $(b_1,\dots,b_n)$ in $\mathcal S$ is $\in R$,
when $\mathcal S$ satisfies the acc,
$(b_1,\dots,b_n)$ in $\mathcal S$ is in $\Gamma$ for
some $\Gamma$ under the dcc, and with each
$(b_1,\dots,b_n)$ in $\mathcal S$ some abridged
type $(b_1',\dots,b_{n'}')$ with bounded $n'$ is in $\mathcal S$.
{\em Abridged\/} means that both types have the same components $b_i\not=0$.
\end{lem}

\begin{proof}
First suppose that each $(b_1,\dots,b_n)$ in $\mathcal S$ is $\in R$ for some fixed finite set of real
numbers $R$. If $\mathcal S$ does not satisfy the acc, then we can find a strictly increasing set of elements $\beta_1,
\beta_2,\dots$ in $\mathcal S$. We can assume that they all have the same size, that is, there is $n$ such that
$\beta_l=(b_1^l,\dots,b_n^l)$. Since $R$ is finite, there are only finitely many such types, a
contradiction.

Now suppose that we have $\mathcal S$ satisfying the acc and other assumptions of the lemma. Let $R\subset
\Gamma$ be the set of all real numbers appearing as a component in some type in $\mathcal S$. It is enough to
prove that $R$ is finite. If $R$ is not finite, then there is a strictly increasing sequence
$\{r_l\}_{l\in}\subset R$ and an infinite set of types $\beta_1, \beta_2,\dots$ in $\mathcal S$ such that $r_l$ is
a component of $\beta_l$. Replacing each $\beta_l$ by an abridged one
$\beta_l=({b_1^l},\dots,{b_{n'}^l})$, we can assume that ${b_1^l}=r_l$. If $n'=1$, then we get
a contradiction. Otherwise, consider types $\lambda_l=({b_2^l},\dots,{b_{n'}^l})$ and use induction on
size and the dcc property of $\Gamma$
to get an infinite increasing subsequence of $\beta_i$.
By construction it is strictly increasing.
This is a contradiction,
because the set $\{\beta_l\}_{l\in\mathbb{N}}$ does not satisfy acc.

$\Box$
\end{proof}

\begin{proof}(of Addendum \ref{add-main-th2})

By Theorem \ref{main-th2}, each set satisfies the acc. Now Lemma \ref{cc} gurantees the existence of $\Gamma_f\subseteq \Gamma$. $\Box$

\end{proof}

We also proved the following.

\begin{cor}\label{4->6}
Under the assumptions of Theorem \ref{main-th2},
(ii) implies (iii), (iv) and (v).
\end{cor}

\begin{proof}
 See the proof of Theorem \ref{main-th2} above.
\end{proof}

\begin{cor} In Theorems \ref{main-th} and \ref{main-th2} and Addenda for $d\leq 4$, we can drop
the assumption on $\LMMP$ in dimension $\leq 4$. \end{cor}

\begin{proof} The LMMP in dimension $4$ follows from ACC of mld's in dimension $4$
[\ref{mld's}, Cor 5].

\end{proof}

\begin{cor}
In (ii-v) of Theorems \ref{main-th} and \ref{main-th2} and Addenda  for $d=3$, we
can drop the assumption on $\LMMP$ in dimension $\leq 3$ and acc for $\mathcal
S_3(\mbox{Mori-Fano cn})$. \end{cor}

\begin{proof}
LMMP and boundedness of Mori-Fano cn varieties are known in dimension $3$ [\ref{log models}][\ref{KMMT}].
\end{proof}

\begin{cor} In (i-v) of Theorems \ref{main-th} and \ref{main-th2} and Addenda for $d=2$, we
can drop the assumption on $\LMMP$ in dimension $\leq 2$, acc for $\mathcal
S_2(\mbox{Mori-Fano cn})$ and $\ACC$ for mld's in dimension $2$. \end{cor}

\begin{proof} The same plus the fact that  $\ACC$ for mld's is known in dimension $2$ [\ref{acc-cod2}][\ref{A}].
\end{proof}


\section{Log twist}

In this section, we introduce a construction which is crucial for us and
which generalizes [reminds] Sarkisov links of Type I and II
[\ref{Matsuki}, Theorem 13-1-1], and we
establish its basic properties.

 \begin{constr}[Log Twist]\label{construction}
Let $X$ be a $d$-dimensional Mori-Fano variety, and $B$ be a boundary
such that $(X,B)$ is klt and noncanonical in codimension $\ge 2$ (noncn for short), and
$K+B\equiv 0$.
Fix a prime b-divisor (exceptional divisor) $E$ such that
$a:=1-e:=\ldis(X,B)=a(E,X,B)$.  Assume the LMMP in dimension $d$. Then there exists (and is unique
for the fixed $E$) the following transformation of $X$ which we call a \emph{log twist}:

$$
\xymatrix{
Y=Y_1\ar[d]^f &\dashrightarrow &Y_2 &\dashrightarrow &\dots &
\dashrightarrow & Y'=Y_n\ar[d]^{f'} \\
X & & & && & X'
}
$$
where $f\colon Y=Y_1\to X$ is an extremal divisorial extraction
of $E$,
all horizontal modifications $Y_i \dashrightarrow Y_{i+1}, i=1,\dots,n-1$,
are extremal  $-E$-flips, and
$f'\colon Y'=Y_n \to X'$ is either a Mori-Fano fibration
with $\dim X'\ge 1$ and the crepant boundary $B_{Y'}$
such that

(1)
$(Y',B_{Y'})$ is klt, and $K_{Y'}+B_{Y'}\equiv 0$, or

an extremal divisorial
contraction of a divisor $E'$ onto
a Mori-Fano variety $X'$ with the crepant boundary $B_{X'}$
such that

(2)
$(X',B_{X'})$ is klt, and $K_{X'}+B_{X'}\equiv 0$.

In addition, the following two facts hold:

(3) If $D$ is an effective divisor on $Y$ which
is (numerically) seminegative over $X$ then its birational
transform $D'$ on $Y'$ is
semipositive over $X'$, and strictly positive
when $D\not=0$.

(4) Thus, if in (3) $D'$ is also seminagative over $X'$,
then $D=D'=0$.
\end{constr}

\begin{warn} Since we consider log discrepancies of the
pair $(X,B)$, the first blowup $Y\to X$ can be in
a terminal and even nonsingular point of $X$.
\end{warn}

\begin{rem-defn}\label{types} We say that the twist
has Type I if  $Y'\to X$ is a fibration.
Otherwise the twist has Type II (cf. [\ref{Matsuki}, Theorem 13-1-1]).
\end{rem-defn}

\begin{lem}[cf. {[\ref{log models}, Theorem 3.1]}]\label{extraction}
Assume the LMMP in dimension $d$.
Let $(X,B)$ be a log pair, and
$E$ be an exceptional prime divisor of $X$ such that

a) (X,B) is klt;

b) $\dim X=d$; and

c) $a=a(E,X,B)<1$.

Then there exists an extraction of $E$:
that is, a contraction
$$
f\colon Y\to X
$$
with the only exceptional divisor $E$ such that $-E$ is ample over $X$.

Moreover, $Y$ and $f$ are unique, and if $X$ is $\mathbb{Q}$-factorial, then so does $Y$ and $f$ is extremal, that is $\rho(Y/X)=1$.
\end{lem}

\begin{proof}
 Let $g\colon W\to X$ be a log resolution of $(X,B)$ such that $E$ is a divisor on $W$. Let $B_W=B^\sim+\sum_{E_i\neq E} E_i$ where $B^\sim$ is the birational transform of $B$ and $E_i$ are the exceptional/$X$ divisors on $W$. Run the LMMP/$X$ on $K_W+B_W$. At the end, we get a model $\overline{W}$ where $K_{\overline{W}}+B_{\overline{W}}$, the pushdown of $K_W+B_W$, is nef (and big)/$X$. By construction, $\overline{W}$ is $\mathbb{Q}$-factorial and $K_{\overline{W}}+B_{\overline{W}}$ is dlt. All the
$E_i$ are contracted/$\overline{W}$
except $E$, by the negativity lemma [\ref{log flips}] which in turn implies that $K_{\overline{W}}+B_{\overline{W}}$ is klt.

In fact, $K_{\overline{W}}+B_{\overline{W}}\equiv -eE/X$ where $e=1-a>0$, and so $-E$ is nef and big/$X$
by construction. Moreover, $K_{\overline{W}}+B_{\overline{W}}$ is semiample/$X$ [\ref{F}, Theorem 7.1], so
$\overline{W}\to X$ can be factored through contractions $h\colon \overline{W}\to Y$ and $f\colon Y\to X$
such that $K_{\overline{W}}+B_{\overline{W}}\equiv 0/Y$ and $K_Y+B_Y$, the pushdown of $K_{\overline{W}}+B_{\overline{W}}$,
is ample/$X$. Since $E\not\equiv 0/X$,
$E$ is a divisor on $Y$ and it is numerically negative/$X$ which
implies that $E$ is the only exceptional/$X$ divisor on $Y$. Note that, $Y$ is the log canonical model
of $K_W+B_W$ [\ref{log models}, Definition 2.1]. This implies the uniqueness of $f$ and $Y$.

Now suppose that $X$ is $\mathbb{Q}$-factorial. Let $D$ be a Weil divisor on $Y$.
If $D=E$, then $D$ is $\mathbb{Q}$-Cartier by construction. If $D\neq E$, let $D'=f_*D$
and $D''=f^*D'$. Since $D'$ is $\mathbb{Q}$-Cartier, so is $D''$. On the other hand $D''=D+\alpha E$
for some rational number $\alpha$. Since $E$ is $\mathbb{Q}$-Cartier, so is $D$. This implies that $Y$
is also  $\mathbb{Q}$-factorial. Moreover, this observationa also shows that $\rho(Y/X)=1$.

$\Box$
\end{proof}

\begin{rem}
We expect that Lemma \ref{extraction} holds in the lc case as well. One then needs not only the LMMP but also the semiampleness for dlt pairs in the birational case.
\end{rem}

\begin{proof}(of Construction \ref{construction})
 There is an extremal extraction $f\colon Y=Y_1\to X$ of $E$
by Lemma \ref{extraction}.
Let $B_Y$ be the crepant pullback of $B$ on $Y$.
According to our assumptions, $B_Y$ is a boundary.
Moreover,

(5) $a<1$ and $0< e\le 1$; and

(6) $\ldis(Y,B_Y)\ge \ldis (X,B)=a$.

By construction $K_Y+B_Y\equiv 0$, and $Y$ is $\mathbb{Q}$-factorial
 by Lemma \ref{extraction}.

Now we run the $\LMMP$ starting from $Y$ with respect
to $K_Y+B_Y-eE\equiv -eE$.
By (5) this is the same as with respect to $-E$.
Since $E$ is always positive on the generic member of some covering family
of curves,
after finitely many flips $Y_i\dashrightarrow Y_{i+1}$,
we get an extremal contraction $f'\colon Y'=Y_{n}\to X'$
which is not a flipping, that is, $f'$ is a Mori-Fano
fibration or
a divisorial contraction, contracting $E'$.
The first case gives a twist of Type I, and
the second one gives a Type II twist.

In both cases,  $E$ is
positive with respect to $f'$, and also so does
$E$ with respect to the flipping contraction
of each flip $Y_i\dashrightarrow Y_{i+1}$.
In particular, $E$ is a divisor on
$X'$ if $f'$ has Type II.
In both cases, the flips are log flops with
respect to $K_{Y_i}+B_{Y_i}$, and all $B_{Y_i}$ are (crepant)
boundaries.
Thus, both Type I and Type II twists satisfy property (1), and in addition,
 the Type II also satisfies (2).
By (6) in both cases,

(6') $\ldis(Y',B_{Y'})\ge \ldis (X,B)=a$.

However, after the contraction in Type II,
this may fail (see Definition \ref{final}).

By construction, $\rho(X')=\rho(Y')-1=
\rho(Y_i)-1=\rho(Y)-1=\rho(X)=1$,
and, for Type II, $X'$ is $\mathbb{Q}$-factorial. 
Hence, for this case,
since $E$ is not exceptional on $X'$ and by (5),
$-K_{X'}$ is ample which means that $X'$
is a Mori-Fano variety.

Now let $D$ be an effective divisor on $Y$ which is seminegative/$X$.
According to the previous paragraph, each $\rho(Y_i)=2$.
Let $R_1$ be an extremal ray corresponding to the contraction
$Y\to X$, and $R_2$ be the other extremal ray.
By our assumption, $D$ is seminegative on $R_1$.
Since $D\ge 0$, $D$ is semipositive on $R_2$.
Thus the first flip $Y_1\dashrightarrow Y_2$
is  a $-D$-flip or $-D$-flop.
Similarly, each next flip $Y_i\dashrightarrow Y_{i+1}$
is a $-D_i$-flip or $-D_i$-flop where $D_i$ denotes
the birational transform of $D$.
This flip corresponds to the second extremal ray $R_2$,
whereas $R_1$ is flipped on $Y_{i+1}$.
So, $D_i$ is always seminegative on $R_1$ and
semipositive on $R_2$.
For $Y'=Y_n$ and $D':=D_n$, this gives the semipositivity in (3) because $R_2$
determines the last contraction $Y'\to X'$.
(This proves also the uniqueness of the twist for
the fixed $E$.)

Moreover, if $D'$ is also seminegative over $X'$, then $D'$ is seminegative on $Y'$ and thus it is $0$
because $D'\ge 0$. Thus $D=0$ too. This proves (4) and the strict positivity in (3). $\Box$ \end{proof}

\begin{defn}\label{final}
A log twist is called \emph{final}, if

(a) $Y'\to X'$ is a fibration, that is, it is of Type I; or

(b) $Y'\to X'$ is of Type II, $X'$ is noncn, and
$\ldis(X',B_{X'})=1-e'$ where $e'=\mult_{E'}B_{Y'}$; or

(c) $Y'\to X'$ is of Type II, and $X'$ is canonical (cn for short).
\end{defn}

Indeed, if a log twist is not final, it is of Type II
with noncn $X'$.
Thus we can take a log twist of $(X',B_{X'})$.
In case (b) of Definition \ref{final}, an inverse log twist
can be constructed. Otherwise we expect that a sequence of
log twists:
$$
(7) ~~~~~~~~~~~~~~(X,B)\dashrightarrow (X',B_{X'})\dashrightarrow  \dots
\dashrightarrow (X^{(i)},B_{X^{(i)}})\dashrightarrow \dots
$$
terminates, where each log twist is nonfinal, except possibly for the last one.

\begin{prop}[Termination of log twists]\label{term-twist}
Suppose that for a sequence as in (7), there exists a real number $a_0<1$ such that

$$\mathrm{(UBD)}~~~~~~~~~~~~~~ each~~ a^{(i)}=mld(X^{(i)},B_{X^{(i)}})\le a_0$$

Then, assuming $\LMMP$ in dimension $d$, the sequence terminates and {\em universally\/}
with respect to $a_0$, that is, the sequence is
finite and the number of twists in it is bounded
whereas the bound depends only on $a_0$ and
the dimension of $X$.
\end{prop}

By $\ACC$ for mld's in dimension $d$ near $1$, we mean that $1$ is not an upper limit in the mld spectrum
(\ref{acc}) in dimension $d$. This is a special case of Conjecture \ref{acc}.

\begin{lem}\label{cor-term} $\ACC$ for mld's in dimension $d$ for
$\Gamma=\{0\}$ and only near $1$
implies \emph{UBD}.

More precisely, there exists $a_0<1$ depending
only on $d$ such that, for any log pair $(X,B)$, with
$\mathbb{Q}$-factorial noncn $X$ of dimension $d$,
$$
a=\ldis(X,B)\le a_0.
$$
\end{lem}

\begin{proof} Put
$$
a_0=\max\{\ldis(X,0)|~ \dim X=d \}\cap [0,1).
$$
Then by ACC for mld's,
$a_0<1$, and for any log pair in the lemma,
$a=\ldis(X,B)\le \ldis(X,0)\le a_0$.
$\Box$
\end{proof}

\begin{cor}\label{cor-term2} Assume  that $\LMMP$ and $\ACC$ for mld's near $1$ hold in dimension
$d$. Let $\Gamma\subset [0,1]$ be a set of real numbers satisfying dcc. Then, universal termination of
Proposition \ref{term-twist} holds in dimension $d$ for $B\in \Gamma$  without the \emph{UBD} assumption,
that is, the length of a sequence of nonfinal twists for $(X,B)$ is bounded by a natural number $N$, where
$\dim X=d$ and $B\in \Gamma$; and $N$ depends only on $d$. \end{cor}

\begin{proof} Immediate by Proposition \ref{term-twist} and Lemma \ref{cor-term}
because each $X$ is $\mathbb{Q}$-factorial.
$\Box$
\end{proof}

\begin{adden}\label{add1} Let
$$\Gamma'=\Gamma\cup\{1-\ldis(X,B)| (X,B)~\mbox{as in (7), $B\in \Gamma$, $X$ noncn and $\dim X=d$}\}$$

{\flushleft and} let $\Gamma^{(i)}=(\Gamma^{(i-1)})'$. Then, the increasing
sequence
$$
\Gamma\subseteq\Gamma'\subseteq\dots\subseteq\Gamma^{(i)}\subseteq\dots
$$
stabilizes, and satisfies the dcc, that is,
the union
$$
\Gamma^{\infty}=\bigcup_1^N\Gamma^{(i)}=\Gamma^{(N)}
$$
and it satisfies the dcc.
\end{adden}

\begin{proof}
More generally, in the addendum,
we can take lc pairs $(X,B)$
with noncn (even non terminal) $X$,
$\dim X=d$, and $B\in \Gamma$ but
we need to assume that $i\le N$, that is,
we state that $\Gamma^{(N)}$ satisfies the dcc.
The assumption holds for the pairs in Proposition \ref{term-twist}
by Corollary \ref{cor-term2}.

Then the new $\Gamma'$ satisfies the dcc as
a union of two sets under the dcc.
The second set
$$\{1-\ldis(X,B)|
(X,B) ~~\text{is lc, but noncn, } \Gamma \in B
\text{ and } \dim X=d\},
$$
satisfies the dcc by $\ACC$ for mld's.
Etc.
$\Box$
\end{proof}

\begin{thm}\label{bnd-support}
 Let $(X/Z,B=\sum b_iD_i)$ be a log pair of dimension $d$ such that:

a) $X\to Z$ is a proper contraction;

b) $(X,B)$ is lc; and

c) $K+B$ is seminegative$/Z$.

Then LMMP in dimension $d$ implies:
$$
\rho^W(X/Z)\ge -\dim X+\sum b_i
$$
where $\rho^W$ is the Weil number, that is,
the rank of Weil divisors on $X$ modulo
numerical equivalence.

\end{thm}

\begin{proof} See [\ref{toric}, Theorem 2.3].
\end{proof}

\begin{cor}\label{cor-bnd-support} Under the assumptions of Theorem \ref{bnd-support}, $\sum b_i\le \dim X +1$ when
$\rho^W(X/Z)=1$.
\end{cor}
Note that $\rho^W=\rho$ for $\mathbb{Q}$-factorial $X$.

\begin{proof} Obvious by Theorem \ref{bnd-support}.
\end{proof}

\begin{proof}(of Proposition \ref{term-twist})
Note that if a log twist is not final then
(6') implies
$$\mbox{(6'')} ~~~~~~~~~~ a'=\ldis(X',B_{X'})=\ldis(Y',B_{Y'})\ge \ldis (X,B)=a$$
and

(8) for the prime divisor $E'$ on $X$,
$1\ge a(E',X',B_{X'})=a(E',X,B)>a'$.

Similarly, for any nonfinal twist
$X^{(i)}\dashrightarrow X^{(i+1)}$,

(8') the prime divisor $E^{(i+1)}$

contracted by $X^{(i)}\dashrightarrow X^{(i+1)}$
satisfies
$$
1\ge a(E^{(i+1)},X^{(i)},B_{X^{(i)}})=
a(E^{(i+1)},X^{(i+1)},B_{X^{(i+1)}})
>a^{(i+1)}.
$$

Thus the sequence of mld's $a,a',\dots, a^{(i)},\dots$
increasing:
$$
a\le a'\le \dots \le a^{(i)}\le \dots,
$$
or equivalently,
$$
(9)~~~~~~e=1-a\ge e'=1-a'\ge \dots\ge e^{(i)}=1-a^{(i)}\ge \dots.
$$

Since the set of boundary multiplicities can be
expanded by twists, it might not be
expected that the sequence $a^{(i)}$ stabilizes or that it is finite.

To establish this we introduce the {\em difficulty\/} $d^{(i)}$ of
the $(X^{(i)},B_{X^{(i)}})$ to be:
the number of prime components $D_i$ of $B_{X^{(i)}}$ with
$b_i=mult_{D_i} B_{X^{(i)}}\ge e^{(i)}$.

The difficulty increases: for any nonfinal twist
$X^{(i)}\dashrightarrow X^{(i+1)}$,
$$
(10)~~~~~ d^{(i+1)}\ge d^{(i)}+1.
$$
Indeed, $e^{(i+1)}\le e^{(i)}$ by (9), and by (6')
none of the prime boundary components $D_i$
of $B_{X^{(i)}}$ with
$b_i=1-a(D_i,X^{(i)},B_{X^{(i)}})\ge e^{(i)}\ge
e^{(i+1)}=1-a^{(i+1)}$
is contracted.
On the other hand, the nonfinal twist of $X^{(i)}$
blows up a new prime component with
the multiplicity $e^{(i)}\ge e^{(i+1)}$
which adds $1$ in the inequality (10).

Note now, that by Theorem \ref{bnd-support} there exists
a natural number $N$ such that,
on each Mori-Fano variety of dimension $d$, a boundary $B$ has
at most
$N$ boundary components with multiplicities $\ge e_0=1-a_0>0$ if $K+B$ is seminegative, and $(X,B)$ is lc. More precisely, we can take any $N\geq (d+1)/e_0$.

By UBD each $a^{(i)}\le a_0$, except for
the last one (if such exists), or equivalently,
each $e^{(i)}\ge e_0=1-a_0$, except for
the last one. Thus we have at most $N$ nonfinal twists.
$\Box$
\end{proof}

\begin{lem}\label{inf-mult} Let $\Gamma\subset [0,1]$ be a set of real numbers satisfying the dcc, and
$X_i\dashrightarrow X_i'$ be a family of birational,
(e.g., nonfinal) log twists in dimension $d$
such that

{\rm a)}  $B_i\in \Gamma$, and

{\rm b)}  the set of multiplicities of all boundaries $B_i$ is infinite.

Then $\ACC$ for mld's and lc thresholds in dimension $d$ imply that
the set of multiplicities of all boundaries $B_{X_i'}$ is infinite too.
\end{lem}

\begin{proof}
By our assumptions each $X_i'$ is a birational modification of $X_i$
with a divisorial contraction $Y_i'\to X_i'$.
In particular, all crepant boundaries $B_{Y_i'}$ and $B_{X_i'}$ are well defined.

By the dcc of $\Gamma$ and after taking a subsequence, we can suppose that
there exists a sequence of prime divisors $D_{i}$ on $X_i$ such
that the corresponding sequence of boundary multiplicities
$b_i=\mult_{D_i} B_i$ is strictly increasing, and the increasing
holds for other multiplicities.
If infinitely many members of the sequence $D_i$ are
nonexceptional on $X_i'$, the required statement holds.

Otherwise, after taking a subsequence, we can suppose that
each $D_i$ is contracted on $X_i'$, that is, $D_i=E_i'$ on $Y_i'$, and
it is numerically negative on $Y_i'$ over $X_i'$.
Thus by the property (3) of twists,
$D_i$ is numerically positive on $Y_i$ over $X_i$.
Hence each $D_i$ passes through $P_i=C_{X_i} E_i$, the center of $E_i$ on $X_i$.
According to $\ACC$ for mld's and lc thresholds, Proposition \ref{mld-th} and the monotonicity of multiplicities,
the set of new boundary multiplicities $e_i=1-\mld(P_i,X_i,B_i)$ is
not finite. But each component $E_i$ is nonexceptional on $X_i'$ which gives the required infinity in
this case too. $\Box$ \end{proof}

\begin{adden}
We can omit $\ACC$ for lc thresholds in Lemma \ref{inf-mult} if we assume the $\LMMP$ and Conjecture \ref{weak-BAB} in dimension $d-1$.
\end{adden}

\begin{proof}
Will be given in section 5.
\end{proof}

\begin{cor}\label{inf-mult1} Assume $\ACC$ for mld's and lc thresholds in dimension $d$. Under the
assumptions of Corollary \ref{cor-term2}   let $(X_i,B_i)$ be a family of pairs in dimension $d$ such
that

{\rm a)}  each has at least one twist,

{\rm b)}  $B_i\in \Gamma$, and

{\rm c)}  the set of multiplicities of all boundaries $B_i$ is infinite.

Then, for their final twists $X_i^{(j)}\dashrightarrow X_i^{(j+1)}$,
the set of multiplicities of all boundaries $B_{X_i^{(j)}}$ is infinite too.
(Taking a subsequence we can suppose also that $j$ is the same
for all final twists.)
\end{cor}

Note that the final twists exist by Corollary \ref{cor-term2}.

\begin{proof} Immediate by Lemma \ref{inf-mult} and induction on $j$. (The statement about a subsequence
follows from the universal termination in Corollary \ref{cor-term2}.) \end{proof}

\begin{adden}
We can omit $\ACC$ for lc thresholds in Corollary \ref{inf-mult1} if we assume the $\LMMP$ and Conjecture \ref{weak-BAB} in dimension $d-1$.
\end{adden}

\begin{proof}
Will be given in section 5.
\end{proof}


\section{Weak finiteness}

  Let $\mathcal{F}_d(\Gamma)$
be the set of log
pairs  $(X,B)$ where $X$ is a $d$-dimensional projective variety, and
$B$ is a boundary such that $B\in \Gamma$,
$(X,B)$ is lc, and $K+B\equiv 0$.

\begin{prop}[Weak finiteness]\label{Weak finiteness}
We assume $\LMMP$, $\ACC$ for mld's and Boundedness
Conjecture \ref{weak-BAB} in dimension $\leq d$.
Let $\Gamma\subset [0,1]$ be a set satisfying
the dcc.
Then there exists a finite subset $\Gamma_f
\subset\Gamma$ such that
$$
\mathcal{F}_d(\Gamma)=\mathcal{F}_d(\Gamma_f),
$$
that is, for each pair $(X,B)\in \mathcal{F}_d(\Gamma)$,
actually $B\in \Gamma_f$.
\end{prop}

\begin{exa}\label{d=1}
Let $\Gamma$ be a set satisfying the dcc,
and $(X,B)\in \mathcal{F}_{1}(\Gamma)$.
So, $X$ is a nonsingular curve, and
since
$K+B\equiv 0$, either $X\iso \mathbb{P}^1$, or
$X$ is an elliptic curve.
In the latter case, $B=0$ and $\Gamma_f=\{0\}$ is enough.

In the
former case, $\deg K+B=0$, and $\sum_1^k b_j=-2$ where
$B=\sum b_j P_j$ and each $b_j\in \Gamma$. By the dcc of $\Gamma$, we may assume that $k$ is bounded, that
is, it only depends on $\Gamma$. Since $(X,B)$ is lc, $k\geq 2$. Note that if $\{s_i\}$ and $\{s'_i\}$ are two sequences satisfying the
dcc, then $\{s_i+s'_i\}$ is also a sequence satisfying the dcc. Similarly, the sum of $n$ dcc sequences,
satisfies the dcc.

Now, if there is no finite $\Gamma_f$ as in Proposition \ref{Weak finiteness}, then there is a sequence
$(X_i,B_i)\in \mathcal{F}_1(\Gamma)$ such that the set of multiplicities
of all $B_i=\sum_1^k b_{i,j} P_{i,j}$
is not finite. In particular, we may assume that all $X_i$ are the projective line and that $\{b_{i,1}\}$
is an infinite dcc sequence. On the hand, $2-b_{i,1}=\sum_2^k b_{i,j}$. Since each $\{b_{i,j}\}$ is a dcc
sequence,  $\sum_2^k b_{i,j}$ is also a dcc sequence. This contradicts the fact that $\{2-b_{i,1}\}$ is an
infinite acc sequence.
\end{exa}

\begin{proof}(of Proposition \ref{Weak finiteness} )

We use induction on $d$.

Step 1. For $d=1$, see Example \ref{d=1}. Now we assume that the theorem holds
in dimension $\le d-1$, and we
establish it in dimension $d$.
Suppose that there exists
a sequence of log pairs $(X_i,B_i)\in F_d(\Gamma)$, $i=1,2,\dots,$
such that the set of boundary multiplicities
$M=\{b_{i,k}\}$, for boundaries $B_i=\sum b_{i,k} D_{i,k}$
is infinite.
Since $M$ satisfies the dcc we can
assume that the sequence $b_{i,1}, i=1,2, \dots$
is strictly increasing, and has only positive
real numbers. Below we derive a contradiction (see Step 8).

Step 2. {\em We can suppose that each $X_i$ is
$\mathbb{Q}$-factorial, and, in particular, each
divisor $D_{i,1}$ is $\mathbb{Q}$-Cartier.\/} Indeed a $\mathbb{Q}$-factorialization
$Y_i\to X_i$ exists by $\LMMP$
[\ref{Isk-Sh}, Corollary~6.7], if $(X_i,B_i)$ is klt,
or the same construction gives
a crepant strictly lt model, and we can replace
each pair $(X_i,B_i)$ by its
crepant $\mathbb{Q}$-factorialization $(Y_i,B_i)$ or
respectively by a crepant strictly log terminal model
$(Y_i,B_i^{\log})$ where the boundary $B_i$ on $Y_i$ is
the birational transform of $B_i$ on $X_i$ or
respectively $B_i^{\log}$ is the log birational transform
of $B_i$. In the last case we need to extend $\Gamma$ by $1$.

Step 3. {\em We can suppose that
each $X_i$ is a Mori-Fano variety\/},
that is, $X_i$ is a projective $\mathbb{Q}$-factorial
variety with Picard number $\rho(X_i)=1$
having only lc singularities and ample $-K_{X_i}$
(cf. [\ref{Isk-Sh}, Definition~1.6 (v)]).

Indeed, by our assumptions we can apply
$\LMMP$ to $(X_i,B_i'=B_i-b_{i,1}D_{i,1})$ or,
equivalently, with respect to the log canonical
divisor  $K_{X_i}+B_i-b_{i,1}D_{i,1}$.
Note that $K_{X_i}+B_i-b_{i,1}D_{i,1}\equiv -b_{i,1}D_{i,1}$
where $b_{i,1}>0$.
Thus $D_{i,1}$ is positive on each extremal
contraction $X_i\to Z_i$.
In particular, the divisor $D_{i,1}$ will never be
contracted.
Moreover, if the extremal contraction is
birational, its birational modification $(X_i^+,(B_i')^+)$
(a divisorial contraction or a log flip;
both are log flips in the sense of [\ref{log models}]) belongs
again to $F_d(\Gamma)$ and we can replace
pair $(X_i,B_i)$ by its log flop $(X_i^+,B_i^+)$.
Of course, the entire set of multiplicities
can decrease but its monotonic infinite subset $\{b_{i,1}\}$
will remain.
On the other hand, we always have an extremal
contraction: $-D_{i,1}$ is not nef, for any $i$.
Therefore, after finitely many steps, the extremal
contraction will be
a Mori-Fano fibration $X_i\to Z_i$.
By construction each $X_i$ has such a contraction.

If $\dim Z_i=d_i\geq 1$, then the generic fibre of
$(X_i/Z_i,B_i)$ with induced (intersection)
boundary belongs to $F_{d-d_i}(\Gamma)$, and
$D_{i,1}$ gives a boundary component with
the same multiplicity $b_{i,1}$.
Hence, by induction we do not have
infinite subsequence $(X_i,B_i)$ with
$d-d_i<d$.
Thus, replacing with a subsequence, we can assume that
each $d_i=0, Z_i=\pt$, and $X_i$ is
a Mori-Fano variety.

Step 4. {\em We can suppose that only finitely
many varieties $X_i$ are cn,\/} and thus,
replacing by a subsequence,
{\em we can suppose that all varieties $X_i$ are noncn.}
Otherwise, we can suppose that each $X_i$ is cn.
Then by Conjecture \ref{weak-BAB}, varieties $X_i$
belong to a bounded family.
Hence, by [\ref{MP}, Lemma 8.1] the set of boundary multiplicities,
in particular,
$\{b_{i,1}\}$ is finite. This is a contradiction.

So, replacing by a subsequence,
we can suppose that each $X_i$ has
a noncn point $P_i$ (it may be nonclosed) and
its codimension $\ge 2$.
We can assume also that each $(X_i,B_i)$ is
klt by Theorem \ref{main-th2} (iii).
Fix a prime b-divisor $E_i$ with center $C_{X_i} E_i=P_i$, and
$\ldis(X_i,B_i)=\mld(P_i,X_i,B_i)=a(E_i,X,B_i)$.
 Thus, we can apply to each $X_i$ a log twist as in Construction \ref{construction}.

Step 5. {\em We can suppose that each log twist
$X_i\dashrightarrow X_i'$ is final.}
Indeed, if it is not final, put $(X_i,B_i)=(X_i',B_{X_i'})$, and
take another twist, etc.
According to our assumptions (ACC of mld's) and Corollary \ref{cor-term2},
after a bounded number of twists, we can suppose that
each twist $X_i\dashrightarrow X_i'$ is final.
By Addendum \ref{add1}, the extended set of boundary multiplicities
$\Gamma^\infty=\Gamma^{(N)}$ again satisfies the dcc.

For the final twist, we denote the resulting contraction  by
$f'\colon Y_i'\to X_i'$.
By Corollary \ref{inf-mult1} there exists an infinite set of distinct boundary
multiplicities for pairs $(X_i,B_i)$.
As in Step 1 and after taking a subsequence, we
can suppose that there exists a sequence of prime divisors
$D_{i,1}$ on $X_i$ with strictly increasing boundary multiplicities
$b_{i,1}$.

Put $a_i=\mld(P_i,X_i,B_i)=a(E_i,X_i,B_i)$
and codiscrepancy $e_i=1-a_i$.
Note that by ACC for mld's and the dcc for $\Gamma$,
the set of
mld's $\{a_i|i=1,2,\dots\}$ satisfies the acc.
Hence the numbers $e_i$ satisfy the dcc, and
we can suppose (after taking a subsequence)
that numbers $e_i$ form a monotonically increasing sequence.
Thus the crepant boundaries
$B_{Y_i}= B_i+e_iE_i$
on $Y_i$, where divisors $B_i$ on $Y_i$
denote the birational transform of $B_i$,
and their modifications $B_{Y_i'}$ belong to the set $\Gamma\cup\{e_i\}$ which
again satisfies the dcc (cf. Addendum \ref{add1}), and the divisors
itself have two monotonically increasing
multiplicities $b_{i,1}$ and $e_i$.
The former is strictly monotonic.

Step 6. {\em Infinitely many twisted contractions
are divisorial,} that is, of Type II (see \ref{types}).
Otherwise infinitely many
of twisted contractions $Y_i'\to X_i'$
are Mori-Fano fibrations.
Then, replacing by a subsequence, we can assume that all of them
are fibrations with varieties $X'_i$ of
the same dimension $\ge 1$.
Thus,  generic fibres will be of the same
dimension $\le d-1$.

On the other hand, it is impossible by induction
that infinitely many birational transforms of
$D_{i,1}$ on $Y_i'$ are strictly positive over $X_i'$,
because then they intersect the generic fibre
(cf. Step~3).

 Therefore, by
(3)
in Construction \ref{construction},
we can assume that each $D_{i,1}$ is positive over $X_i$.
Then, by Proposition~\ref{mld-th}
the set $\{e_i\}$ is not finite, and
by ACC for mld's
we can suppose that $e_i$ is strictly increasing (cf.
Step~8 below and the proof of Lemma \ref{inf-mult}).
This again gives a contradiction because,
by
(3)
in Construction \ref{construction}, the birational
transform of each $E_i$ on $Y_i'$ is positive over $X_i'$.

Thus we can assume that each twisted contraction
$Y_i'\to X_i'$ is divisorial with $X_i'$ a Mori-Fano variety, and
some divisor $E_i'$ is contracted to a point $P_i'$ in $X_i'$ of
codimension $\ge 2$.
Put $a_i'=a(E_i',X_i',B_{X_i'})$ and
the codiscrepancy $e_i'=1-a_i'$.

Step 7: {\em We can assume that $a_i'=\mld(P_i',X_i',B_{X_i'})$ for infinitely many $i$, so we can assume
that for all $i$}. Otherwise, since the twists are final, $X'_i$ is canonical for infinitely many $i$. This
is a contradiction, because such varieties are bounded.

Step 8. {\em Contradiction: $M$ is finite.}
By the dcc of $\Gamma$ and since the support of
$B_i$ has a bounded number of components (Corollary \ref{cor-bnd-support}), we can assume that each sequence
$b_{i,k}, i=1,2,\dots$, is increasing with respect to $i$.
The crepant divisors
$$
B_{Y_i}=\sum_{k=1}^{n+1}b_{i,k} D_{i,k}=e_i E_i+\sum_{k=1}^n b_{i,k} D_{i,k}
$$
satisfy the same property where $b_{i,n+1}:=e_i$. Now we define the
set $R=\{r_k|k=1,2,\dots,n+1\}$ as the set of limits (not necessarily distinct)
$$
r_k=\lim _{i\to \infty} b_{i,k}
$$
First suppose that $a=\lim_{i\to \infty} a_i>0$ and $a'=\lim_{i\to \infty} a_i'>0$. Let $\tau$ be a positive real number
constructed in Main Proposition \ref{transform}. We can assume that each $b_{i,k}\in [r_k-\tau,r_k]$.
Hence by Proposition \ref{aa} each $K_{Y_i}+B_{Y_i}^\tau$ is seminegative over $X_i$, and so
does its birational transform $K_{Y_i'}+B_{Y_i'}^\tau$ over $X_i'$. Here and for the rest of
the proof, the superscript $\tau$ stands for the limit, that is, for example $B_{Y_i}^\tau=\lim
_{j\to \infty} B_{Y_j}$ on $Y_i$ in the sense of Example \ref{bb}. On the other hand, by
construction
$$
B_{Y_i}^\tau\ge B_{Y_i}
.
$$
Hence

$$
D:=K_{Y_i}+B_{Y_i}^\tau= K_{Y_i}+B_{Y_i}+
(B_{Y_i}^\tau -B_{Y_i})
\equiv B_{Y_i}^\tau -B_{Y_i}\ge 0
$$

{\flushleft is} an effective divisor which is
seminegative over $X_i$.
Since flips preserve numerical equivalence, $D'$
the birational transform of $D$ is
seminegative over $X_i'$.
Thus, by (3) in Construction \ref{construction},
$D=D'=0$.
This means that all limits are stabilized.
This contradicts infinity of $M$, in particular
strict monotonicity of $b_{i,1}$.

Now if $a=0$, then to get the seminegativity of $K_{Y_i}+B_{Y_i}^\tau$ over $X_i$, we can use
Theorem \ref{main-th} (vi). In fact, if $K_{Y_i}+B_{Y_i}^\tau$ is not seminegative over $X_i$,
then $K_{X_i}+A_{i}$ is maximally $0$-lc at some point of $X_i$ for some $B_i\lneq
A_i\lneq B_{i}^\tau$ which contradicts Theorem \ref{main-th2} (vi). We have a similar argument
for $a'=0$. The rest of the proof is exactly as in the $a,a'>0$ case.

$\Box$

\end{proof}

\section{Proof of Main Theorem}

\begin{proof}(of Main Theorem \ref{result})

(i) By induction, we can assume ACC for lc thresholds in dimension $\le d$. Now we can use Proposition \ref{mld-th}.

(ii) This follows from Addendum \ref{main-th-adden} (v).

(iii) We can use (ii) and the main result of [\ref{B}].

(iv) This follows from [\ref{F}] (see [\ref{BCHM}] or [\ref{ordered}] for more information).

$\Box$
\end{proof}

{\noindent DPMMS, Centre for Mathematical Sciences,\\
Wilberforce Road, Cambridge CB3 0WB, UK\\ e-mail:
c.birkar@dpmms.cam.ac.uk}\\

{\noindent Department of Mathematics,\\ Johns Hopkins University,\\
Baltimore, MD--21218, USA\\ e-mail: shokurov@math.jhu.edu}\\

{\noindent
Steklov Mathematical Institute,\\ Russian Academy of Sciences,\\
Gubkina str. 8, 119991, Moscow, Russia}

\end{document}